\documentclass[a4paper,reqno]{amsart}
\usepackage{graphics}
\usepackage{amsmath}
\usepackage{latexsym}
\usepackage{amssymb}
\usepackage[all]{xy}
\xyoption{matrix}
\xyoption{arrow}

\begin{document}

\renewcommand{\mod}{\operatorname{mod}\nolimits}
\newcommand{\fl}{\operatorname{f.l.}\nolimits}
\newcommand{\gl}{\operatorname{gl.dim}\nolimits}
\newcommand{\CM}{\operatorname{CM}\nolimits}
\newcommand{\add}{\operatorname{add}\nolimits}
\newcommand{\Sub}{\operatorname{Sub}\nolimits}
\newcommand{\fd}{\operatorname{f.l.}\nolimits}
\newcommand{\Hom}{\operatorname{Hom}\nolimits}
\newcommand{\RHom}{\operatorname{RHom}\nolimits}
\newcommand{\End}{\operatorname{End}\nolimits}
\renewcommand{\Im}{\operatorname{Im}\nolimits}
\newcommand{\Ker}{\operatorname{Ker}\nolimits}
\newcommand{\Cok}{\operatorname{Cok}\nolimits}
\newcommand{\Ext}{\operatorname{Ext}\nolimits}
\newcommand{\Tor}{\operatorname{Tor}\nolimits}
\newcommand{\op}{{\operatorname{op}}}
\newcommand{\id}{{\operatorname{id}}}
\renewcommand{\k}{K}
\newcommand{\C}{\operatorname{\mathcal C}\nolimits}
\newcommand{\D}{\operatorname{\mathcal D}\nolimits}
\newcommand{\E}{\operatorname{\mathcal E}\nolimits}
\newcommand{\T}{\operatorname{\mathcal T}\nolimits}
\newcommand{\R}{\operatorname{\mathbb R}\nolimits}
\renewcommand{\P}{\operatorname{\mathcal P}\nolimits}
\newcommand{\JJ}{\operatorname{\mathcal J}\nolimits}
\newcommand{\Min}{M_{\operatorname{in}}}
\newcommand{\Mout}{M_{\operatorname{out}}}
\newcommand{\fin}{\operatorname{f_{in}}\nolimits}
\newcommand{\fout}{\operatorname{f_{out}}\nolimits}
\renewcommand{\top}{\operatorname{top}\nolimits}
\newcommand{\pd}{{\operatorname{pd}}}
\newcommand{\F}{\operatorname{\mathbb F}\nolimits}

\newtheorem{lemma}{Lemma}[section]
\newtheorem{proposition}[lemma]{Proposition}
\newtheorem{definition}[lemma]{Definition}
\newtheorem{cor}[lemma]{Corollary}
\newtheorem{theorem}[lemma]{Theorem}

\author[Buan]{Aslak Bakke Buan}
\address{Institutt for matematiske fag\\
Norges teknisk-naturvitenskapelige universitet\\
N-7491 Trondheim\\
Norway}
\email{aslakb@math.ntnu.no}

\author[Iyama]{Osamu Iyama}
\address{
Graduate School of Mathematics \\ Nagoya University \\ Chikusa-ku \\  Nagoya \\  464-8602 \\ Japan}
\email{iyama@math.nagoya-u.ac.jp}

\author[Reiten]{Idun Reiten}
\address{Institutt for matematiske fag\\
Norges teknisk-naturvitenskapelige universitet\\
N-7491 Trondheim\\
Norway}
\email{idunr@math.ntnu.no}

\author[Smith]{David Smith}
\address{Department of Mathematics \\ Bishop's University \\  Sherbrooke, QC \\ Canada \\ J1M1Z7}
\email{dsmith@ubishops.ca}

\title{Mutation of cluster-tilting objects and potentials}

\begin{abstract}We prove that mutation of cluster-tilting objects in triangulated 2-Calabi-Yau 
categories is closely connected with mutation of quivers with potentials. This gives a close connection between
2-CY-tilted algebras and Jacobian algebras associated with quivers with potentials.

We show that cluster-tilted algebras are Jacobian and also that they are determined by their quivers. 
There are similar
results when dealing with tilting modules over 3-CY algebras.
The nearly
Morita equivalence for 2-CY-tilted algebras is shown to hold for the finite length modules
over Jacobian algebras.
\end{abstract}


\thanks{All authors were supported by the Storforsk-grant 167130 from the Norwegian Research Council,
the second author also by JSPS Grant-in-Aid for Scientific Research 18740007,
the third author also by a Humboldt Research Prize, and the fourth author also by NSERC of Canada.}
\maketitle

\tableofcontents

\section*{Introduction}

The Fomin-Zelevinsky mutation of quivers plays an important role
in the theory of cluster algebras initiated in \cite{fz1}.
There is, motivated by this theory via \cite{mrz}, a
mutation of cluster-tilting objects in cluster categories, and more generally Hom-finite 
triangulated 2-Calabi-Yau (2-CY for short) categories over an algebraically closed field $K$ \cite{bmrrt,iy}. 
This has turned out to give a categorical
model for the quiver mutation 
in certain cases \cite{bmr2,birs}.
There is also a theory of mutation for 3-CY algebras \cite{ir}.
On the other hand there is the recent theory of mutation of quivers with potentials $(Q,W)$,
initiated in \cite{dwz}.

Associated with cluster-tilting objects $T$ in 2-CY categories
$\C$ are the endomorphism algebras $\End_{\C}(T)$, called
2-CY-tilted algebras. And associated with a quiver
with potential $(Q,W)$ are algebras called Jacobian
algebras in \cite{dwz}. The mutation of cluster-tilting objects 
induces an operation on the associated 2-CY-tilted algebras
and the mutation
of quivers with potentials induces an operation on the associated 
Jacobian algebras.
The theme of this paper is to investigate these classes of algebras and 
their relationships, in particular with respect to mutation.
In addition to cluster-tilting objects in triangulated 2-CY categories,
we also deal with mutation of tilting modules over 3-CY algebras, and their 
relationship to mutation of quivers with potentials. We now state the main results
of this paper, referring to section \ref{sec1} for definitions 
and background material.

(A) Starting with a tilting module $T$ over a 3-CY algebra $\Lambda$
such that $\End_{\Lambda}(T)$ is Jacobian, our first main theorem states that
the two mutations ``coincide'' (Theorem \ref{main2}). The same type of result is proved for 
cluster-tilting objects in triangulated 2-CY categories (Theorems \ref{main} and \ref{main3}).
Our basic idea for the proof is to use the relationship between Jacobian algebras and 2-almost split sequences introduced in \cite{i,iy},
which is a higher analogue of the classical relationship between mesh categories of translation quivers and almost split sequences \cite{R,BG,IT1,IT2,I2}.

(B) The results in (A) can be used to show that large classes  
of 2-CY-tilted algebras are Jacobian. In particular, a main result is that 
cluster-tilted algebras, which are an important subclass of 2-CY-tilted algebras,
are Jacobian (Corollary \ref{jacob}). This is an easy consequence of Theorem \ref{main3}.

More generally we show that a large class of 2-CY-tilted algebras coming from 
triangulated 2-CY categories associated with elements in Coxeter groups (see \cite{birs})
are also Jacobian (Corollary \ref{mutation of reduced words are given by QP}).

(C) It is an open question whether there exists a
(mutation-) operation on algebras with the following property: $\End_{\C}(T)$ 
should be sent to $\End_{\C}(T^{\ast})$, when $T$ is 
a cluster-tilting object in a triangulated 2-CY category and $T^{\ast}$
is a cluster-tilting object obtained from $T$ by mutation. 

It is an easy consequence of (A) that this question has an affirmative answer 
under some additional conditions.

(D) Another main result in this paper is that cluster-tilted algebras
are determined by their quivers. This was shown for finite representation
type in \cite{bmr3}. 
We here give a short proof of this fact (Theorem \ref{thm2.3}). Alternatively, it is
also a direct consequence of (B).

(E) One of the starting points of tilting theory is the reflection functors \cite{bgp, apr}, giving a 
{\em nearly Morita equivalence}, in the terminology of \cite{ringel}, meaning that 
$\frac{\mod \Lambda}{[\add S_k]} \simeq \frac{\mod \Lambda'}{[\add S_k']}$, where $S_k$ and $S_k'$ are the simple modules
associated with the vertex $k$, where $k$ is a sink or a source in the quiver of the path algebra $\Lambda$,
and $\Lambda'$ is the path algebra obtained by changing direction of all arrows adjacent to $k$.
An important property of cluster-tilted algebras
(and more generally 2-CY-tilted algebras), is that this nearly Morita equivalence generalizes
in the sense that one replaces reflection at a source/sink with mutation at any vertex \cite{bmr1,kr1}.
A consequence of this is
that the cluster-tilted algebras (in one mutation class) have the same 
``number'' of indecomposable modules. In view of the close connection between 2-CY-tilted algebras and
Jacobian algebras it is natural to investigate if a property known to hold for one of the classes 
also holds for the other one.
In this spirit we show that nearly Morita equivalence holds also for two Jacobian algebras 
where one is obtained from the other by mutation of quivers with potentials (Theorem \ref{nearly}). 
We show this also for algebras which are not finite dimensional,
dealing with the category of finite length modules. 
In view of our previous results, this gives a generalization of the above result for cluster-tilted algebras.
For the proof we use a functorial approach to results and techniques of \cite{dwz}. 


The organisation of the paper is as follows. In section 1 we collect some known results on the different kinds of mutation
and on cluster-tilted algebras. In section 2 we prove that cluster-tilted algebras are determined by
their quiver. The connection between tilting and cluster-tilting mutation and mutation of quivers with potential is given in sections 3, 4 and 5.
In section 6 we show that a class of 2-CY-tilted algebras associated with reduced expressions 
of elements of the Coxeter groups are given by QP's, which are rigid in the sense of
\cite{dwz}. The result on nearly Morita equivalence is given
in section 7. 

The results in this paper have been presented at conferences in Oxford and Oberwolfach \cite{reiten-o, i3}.
That cluster-tilted algebras and some of the 2-CY-tilted algebras associated with elements in Coxeter groups are 
Jacobian is proved independently in \cite{k3} using completely different methods (see also \cite{amiot}).
The result on nearly Morita equivalence is independently proved in \cite{dwz2}.
Our result (A) on 3-CY algebras is related to results of Vit\'oria \cite{v} and Keller-Yang \cite{ky}.

\medskip
\noindent{\bf Conventions }
All modules are left modules, and a composition $ab$ of morphisms (respectively, arrows) $a$ and $b$ means first $a$ and then $b$.
We denote by $s(a)$ the start vertex of an arrow or path $a$ and $e(a)$ denotes the end vertex.
For an algebra $\Lambda$ we denote by $\mod\Lambda$ the category of finitely generated $\Lambda$-modules, and by $\fl\Lambda$ the category of finite length $\Lambda$-modules.
We denote by $J_\Lambda$ the Jacobson radical of an algebra $\Lambda$, and by $J_{\C}$ the Jacobson radical of an additive category $\C$.
For an object $T$ in an additive category $\C$ we denote by $Q_T$ the quiver of $\End_{\C}(T)$.

\section{Preliminaries on mutation}\label{sec1}

In this section we discuss different kinds of mutations. We first recall the Fomin-Zelevinsky mutation
of quivers and the recent mutation of quivers with potentials from \cite{dwz}.
Then we consider the mutation of cluster-tilting objects in 2-Calabi-Yau (2-CY) categories, which often gives
a categorical modelling of quiver mutation. 
A main theme of this paper is the comparison of the last two mutations. Since
cluster-tilted algebras play a central role in this paper, we recall their definition
and some of their basic properties. Finally, we consider mutation of tilting modules over 3-Calabi-Yau algebras,
which will also be compared to mutation of quivers with potentials.

\subsection{Mutation of quivers}

Let $Q$ be a finite quiver with 
vertices $1, \dots, n$, and having no loops or $2$-cycles. For any vertex $k$, Fomin-Zelevinsky \cite{fz1} defined a new quiver
$\mu_k(Q)$ as follows. Let $b_{ij}$ and $b_{ij}'$ denote the number of arrows from $i$ to $j$ minus the number of arrows from
$j$ to $i$ in $Q$ and $\mu_k(Q)$, respectively. Then we have
$$b_{ij}' = \begin{cases} -b_{ij} & \text{ if } i=k \text{ or } j=k,\\
b_{ij} + \frac{|b_{ik}|b_{kj}+b_{ik}|b_{kj}|}{2} & \text{  else.} \end{cases}$$
Clearly we have $\mu_k(\mu_k(Q))\simeq Q$.
This mutation is an essential ingredient in the theory of cluster algebras initiated in \cite{fz1}.
Note that if the vertex $k$ is a sink or a source, then this mutation coincides with the 
Bernstein-Gelfand-Ponomarev reflection of quivers \cite{bgp}.

\subsection{Mutation of quivers with potentials}
Let $Q$ be a finite connected quiver without loops, and with vertices $1, \dots, n$,
and set of arrows $Q_1$.
We denote by $KQ_i$ the $K$-vector space with basis $Q_i$ consisting of paths of length $i$ in $Q$,
and by $KQ_{i,{\rm cyc}}$ the subspace of $KQ_i$ spanned by all cycles.
Consider the complete path algebra
\[\widehat{KQ}=\prod_{i\ge0}KQ_i\]
over an algebraically closed field $K$.
A \emph{quiver with potential} ({\em QP} for short) is a pair $(Q,W)$
consisting of a quiver $Q$ without loops and an element $W\in\prod_{i\ge2}KQ_{i,{\rm cyc}}$ (called a {\em potential}).
It is called \emph{reduced} if $W\in\prod_{i\ge3}KQ_{i,{\rm cyc}}$.
The cyclic derivative $\partial_aW$ is defined by
$\partial_a(a_1\cdots a_\ell)=\sum_{a_i=a}a_{i+1}\cdots a_\ell a_1\cdots a_{i-1}$ and extended linearly and continuously.
A QP gives rise to what has been called the associated \emph{Jacobian algebra} \cite{dwz}
\[\P(Q,W)=\widehat{\k Q}/\JJ(W),\]
where $\JJ(W)=\overline{\langle \partial_a W \mid a \in Q_1 \rangle}$
is the closure of the ideal generated by $\partial_aW$ with respect to the $J_{\widehat{KQ}}$-adic topology.

Two potentials $W$ and $W'$ are called \emph{cyclically equivalent} if $W-W'\in\overline{[KQ,KQ]}$,
where $[-,-]$ denotes the vector space spanned by commutators.
Two QP's $(Q,W)$ and $(Q',W')$ are called \emph{right-equivalent} if $Q_0=Q'_0$ and there exists a $K$-algebra isomorphism
$\phi \colon \widehat{KQ}\to\widehat{KQ'}$ such that $\phi|_{Q_0}={\rm id}$ and $\phi(W)$ and $W'$ are cyclically equivalent. In this case $\phi$ induces an isomorphism $\P(Q,W)\simeq\P(Q',W')$.
It is shown in \cite{dwz} that for any QP $(Q,W)$ there exists a reduced QP $(Q_{\rm red}, W_{\rm red})$
such that $\P(Q,W) \simeq \P(Q_{\rm red}, W_{\rm red})$, which is uniquely determined up to right-equivalence.
We call $(Q_{\rm red},W_{\rm red})$ a {\em reduced part} of $(Q,W)$.
For example, a reduced part of the QP $(Q,W)$ below is given by the QP $(Q_{\rm red}, W_{\rm red})$ below.
\[(Q,W)=\biggl(\xymatrix@C0.4cm@R0.2cm{{\scriptstyle 1}\ar@/^-1.25pc/[rr]^{c}\ar[r]^{\scriptstyle a}&{\scriptstyle 2}\ar[r]^{\scriptstyle b}&{\scriptstyle 3}\ar@/^-1.25pc/[ll]^{\scriptstyle d}},cd+abd\biggr)\ \ \ \ \ 
(Q_{\rm red}, W_{\rm red})=\biggl(\xymatrix@C0.4cm@R0.2cm{{\scriptstyle 1}\ar[r]^{\scriptstyle a}&{\scriptstyle 2}\ar[r]^{\scriptstyle b}&{\scriptstyle 3}},0\biggr)\]
A QP $(Q,W)$ is called {\em rigid} \cite{dwz} if the deformation space
$\P(Q,W)/ (KQ_0 + \overline{[\P(Q,W),\P(Q,W)]})$
is zero, or equivalently
$J_{\widehat{KQ}} = \JJ(W) + \overline{[KQ,KQ]}$ holds.

A mutation $\mu_k(Q,W)$ of a QP $(Q,W)$ is introduced in \cite{dwz} for each vertex $k$ in $Q$ not lying on a 2-cycle.
It is defined as a reduced part of $(Q',W') = \widetilde{\mu_k}(Q,W)$, the latter one being given as follows.
Replacing $W$ by a cyclically equivalent potential, we assume that no cycles in $W$ start at $k$.
\begin{itemize}
\item[(a)] $Q'$ is a quiver obtained from $Q$ by the following changes.
\begin{itemize}
\item[(i)] Replace the fixed vertex $k$ in $Q$ by a new vertex $k^{\ast}$.
(Although $k$ and $k^{\ast}$ are identified in \cite{dwz}, we distinguish them to avoid any confusion.)
\item[(ii)] Add a new arrow $[ab] \colon i\to j$ for each pair of arrows $a \colon i\to k$ and $b \colon k\to j$ in $Q$.
\item[(iii)] Replace each arrow $a \colon i\to k$ in $Q$ by a new arrow $a^{\ast} \colon k^{\ast}\to i$.
\item[(iv)] Replace each arrow $b \colon k\to j$ in $Q$ by a new arrow $b^{\ast} \colon j\to k^{\ast}$.
\end{itemize}
\item[(b)] $W'=[W]+\Delta$ where $[W]$ and $\Delta$ are the following:
\begin{itemize}
\item[(i)] $[W]$ is obtained by substituting $[ab]$ for each factor $ab$ in $W$ with $a \colon i\to k$ and $b \colon j\to k$.
\item[(ii)] ${\displaystyle\Delta=\sum_{a,b\in Q_1\atop e(a)=k=s(b)}a^{\ast}[ab]b^{\ast}}$.
\end{itemize}
\end{itemize}
Then $k^{\ast}$ is not contained in any 2-cycle in $\mu_k(Q,W)$, and $\mu_{k^{\ast}}(\mu_k(Q,W))$ is right-equivalent to $(Q,W)$ \cite{dwz}.
Note that if both $Q$ and the quiver part of $\mu_k(Q,W)$ have no loops and 2-cycles, then they are in the relationship of Fomin-Zelevinsky mutation.

\medskip
For example, we calculate $\mu_2(Q,W)$ and $\mu_{2^{\ast}}(\mu_2(Q,W))$ for the QP $(Q,W)$ below.
(For simplicity we denote $a^{**}$ and $b^{**}$ by $a$ and $b$ respectively.)
\[\begin{array}{cccc}
&(Q,W)=\biggl(\xymatrix@C0.4cm@R0.2cm{{\scriptstyle 1}\ar[r]^{\scriptstyle a}&{\scriptstyle 2}\ar[r]^{\scriptstyle b}&{\scriptstyle 3}},0\biggr)&\stackrel{\mu_2}{\longrightarrow}&
\biggl(\xymatrix@C0.4cm@R0.2cm{{\scriptstyle 1}\ar@/^-1.25pc/[rr]^{[ab]}&{\scriptstyle 2^*}\ar[l]_{\scriptstyle a^*}&{\scriptstyle 3}\ar[l]_{\scriptstyle b^*}},a^*[ab]b^*\biggr)\\
\stackrel{\widetilde{\mu}_{2^{\ast}}}{\longrightarrow}&
\biggl(\xymatrix@C0.4cm@R0.2cm{{\scriptstyle 1}\ar@/^-1.25pc/[rr]^{[ab]}\ar[r]^{\scriptstyle a}&{\scriptstyle 2}\ar[r]^{\scriptstyle b}&{\scriptstyle 3}\ar@/^-1.25pc/[ll]_{\scriptstyle [b^*a^*]}},[ab][b^*a^*]+b[b^*a^*]a\biggr)&\xrightarrow{\rm reduced}&
\biggl(\xymatrix@C0.4cm@R0.3cm{{\scriptstyle 1}\ar[r]^{\scriptstyle a}&{\scriptstyle 2}\ar[r]^{\scriptstyle b}&{\scriptstyle 3}},0\biggr)
\end{array}\]

We introduce the following generalization of Jacobian algebras of QP's, which play an important role in this paper.

\begin{definition}
We call a triple $(Q,W,F)$ a \emph{QP with frozen vertices} if $(Q,W)$ is a QP and $F$ is a subset of $Q_0$.
We define the associated \emph{Jacobian algebra} by
\[\P(Q,W,F)=\widehat{\k Q}/\JJ(W,F),\]
where $\JJ(W,F)$ is the closure
\[\JJ(W,F)=\overline{\langle \partial_a W \mid a \in Q_1,\ s(a)\notin F\ \mbox{ or }\ e(a)\notin F\rangle}\]
with respect to the $J_{\widehat{KQ}}$-adic topology.
This means that we exclude cyclic derivatives associated with arrows between frozen vertices.
\end{definition}

\subsection{Mutation of cluster-tilting objects}\label{subsec:mut_obj}

Let $\C$ be a $\Hom$-finite triangulated $K$-category, where $K$ is an algebraically closed field. 
We denote by $[1]$ the shift functor in $\C$, and $\Ext^i_{\C}(A,B) = \Hom_{\C}(A,B[i])$.
Then $\C$ is said to be {\em $n$-Calabi-Yau} ({\em $n$-CY} for short)
if there is a functorial isomorphism 
$$D\Hom_{\C}(A,B) \simeq \Ext^n_{\C}(B,A)$$
for $A,B$ in $\C$ and $D = \Hom_K(-,K)$.
Note that this is called weakly $n$-Calabi-Yau in \cite{k2}.
Similarly we call a Frobenius category $\E$ \emph{$n$-CY} if the stable category
$\underline{\E}$ is $n$-CY.

Let $\C$ be 2-CY.
An object $T$ in $\C$ is {\em cluster-tilting} if
\[\add T=\{X\in\C\ |\ \Ext^1_{\C}(T,X)=0\}\]
(see \cite{bmrrt,i, kr1}).
In this case the algebra $\End_{\C}(T)$ is called a {\em 2-CY-tilted} algebra.
We define a \emph{cluster-tilting object} in $\E$ by the same formula.

Let $T = T_1 \oplus T_2 \oplus \cdots \oplus T_n$, where 
the $T_i$ are nonisomorphic indecomposable objects. For each $k= 1, \dots, n$ there 
is a unique indecomposable object $T_k^{\ast}$ in $\C$ with $T_k^{\ast} \not \simeq T_k$ such that
$(T/T_k) \oplus T_k^{\ast}$ is a cluster-tilting object in $\C$, and we write 
$\mu_k(T) = (T/T_k) \oplus T_k^{\ast}$.
Clearly we have $\mu_k(\mu_k(T))\simeq T$.
There are associated {\em exchange triangles}
$$T_k^{\ast} \overset{g}{\to} U_k  \overset{f}{\to} T_k \to T_k^{\ast}[1]\ \mbox{ and }\ 
T_k \overset{g'}{\to} U_k'  \overset{f'}{\to} T_k^{\ast} \to T_k[1]$$
where $f$ and $f'$ are minimal right $\add(T/T_k)$-approximations, and 
$g$ and $g'$ are minimal left $\add(T/T_k)$-approximations \cite{bmrrt,iy}.
Also for a cluster-tilting object in $\E$, we have \emph{exchange sequences} defined similarly.

In general, for a category $\C$ and a full subcategory $\C'$, we denote by $[\C'](X,Y)$ the subgroup
of $\Hom_{\C}(X,Y)$ consisting of morphisms factoring through objects in $\C'$.
Then $[\C']$ forms an ideal of the category $\C$.

We have an equivalence $\Hom_{\C}(T,-)\colon \C/[\add T[1]] \simeq \mod \End_{\C}(T)$ by \cite{bmr1,kr1}.
Using this, if there is no loop at $k$ in $Q_T$, we have an equivalence 
$$\frac{\mod \End_{\C}(T)}{[\add S_k]} \simeq \frac{ \mod \End_{\C}(\mu_k(T))}{[\add S_k^{\ast}]},$$
where $S_k$ and $S_k^{\ast}$ are the simple modules associated with 
$T_k$ and $T_k^{\ast}$.
Thus $\End_{\C}(T)$ and $\End_{\C}(\mu_k(T))$ are by definition nearly Morita equivalent.

We say that a 2-CY category $\C$ has {\em no loops or 2-cycles}
if the quiver $Q_T$ of $\End_{\C}(T)$ has no loops or 2-cycles for any cluster-tilting object $T$ in $\C$ (see \cite{birs}).
Under this assumption we have the following important connection with the Fomin-Zelevinsky quiver mutation.
For a cluster-tilting object $T$ we have
\[Q_{\mu_k(T)}\simeq\mu_k(Q_T)\]
by \cite{bmr2,birs}.
Hence $\C$ has a cluster structure in the sense of \cite{birs}.

\subsection{Cluster-tilted algebras}
Cluster categories are by definition the orbit categories $\C_Q = \D^{\rm b}(KQ)/\tau^{-1}[1]$,
where $Q$ is a finite connected acyclic quiver, $\D^{\rm b}(KQ)$ is the bounded derived category of the finite dimensional
(left) $KQ$-modules, and $\tau$ is the AR-translation in $\D^{\rm b}(KQ)$, see \cite{bmrrt}. These orbit categories
are known to be triangulated by \cite{k1},
and are $\Hom$-finite 2-CY \cite{bmrrt}.
The {\em cluster-tilted} algebras are by definition the 2-CY-tilted algebras coming
from cluster categories \cite{bmr1}. 

For cluster categories $\C_Q$, the cluster-tilting objects are induced by tilting $KQ'$-modules
over path algebras $KQ'$ derived equivalent to $KQ$. It is known \cite{bmrrt} that given any
two cluster-tilting objects in $\C_Q$, then one can be obtained from the other by a finite sequence of mutations
of cluster-tilting objects, and hence there is only one mutation class of cluster-tilted algebras by \cite{bmrrt}, based on \cite{hu}. 

By \cite{abs} (\cite{br} in the Dynkin case) the quiver of $\End_{\C_Q}(T)$ is obtained from the
quiver of the tilted algebra $\End_{KQ'}(T)$ by adding an arrow in the opposite 
direction for each relation in a minimal set of relations. In particular, we have the following, which is also easily seen directly.

\begin{lemma}\label{l:acyc}
Let $T$ be a tilting module in $\mod KQ$. Then $\End_{\C_Q}(T) \simeq KQ'$ for a quiver $Q'$ with no oriented
cycles, if and only if $\End_{KQ}(T) \simeq KQ'$.
\end{lemma}

If the quiver $Q_T$ of a 2-CY-tilted algebra coming from a 2-CY category
$\C$ is acyclic, then it follows from \cite{kr1} that $\End_{\C}(T) \simeq KQ_T$.
Moreover if $\C$
is a connected algebraic triangulated category (that is, the stable category of
a Frobenius category), then
$\C$ is equivalent to the cluster category $\C_{Q_T}$ \cite{kr2}.  

\subsection{Mutation of tilting modules over 3-Calabi-Yau algebras}\label{subsec:mut_obj2}
Let $R$ be a formal power series ring $K[[x,y,z]]$ in three variables over an algebraically closed field $K$,
and let $\Lambda$ be an $R$-algebra which is a finitely generated $R$-module.
We call $\Lambda$ \emph{3-Calabi-Yau} (\emph{3-CY} for short) if the bounded derived category $\D^{\rm b}(\fl\Lambda)$ is a 3-CY category \cite{b,cr,gin,ir}.
It was shown in \cite{ir} that $\Lambda$ is 3-CY if and only if $\Lambda$ is a free $R$-module which is a symmetric 
$R$-algebra with $\gl\Lambda=3$.
Rickard proved that 3-CY algebras are closed under derived equivalences (see \cite{ir}).

Let $\Lambda$ be a 3-CY algebra. For a basic tilting $\Lambda$-module $T$ of projective 
dimension at most one, we have another 3-CY algebra $\Gamma =\End_\Lambda(T)$.
We take an indecomposable decomposition $T=T_1\oplus\cdots\oplus T_n$.
Let $e_k\in\Gamma$ be the idempotent corresponding to $T_k$ for $k=1,\dots,n$.
We let $I_k =\Gamma(1-e_k)\Gamma$.
Let us recall the following result \cite[Th. 5.4, Th. 7.1]{ir}.

\begin{proposition}\label{recall results}
If $\Gamma/I_k\in\fl\Gamma$, then there exists a minimal projective resolution
\begin{equation}\label{projective resolution of Gamma/I}
0\to\Gamma e_k\xrightarrow{f_{2}}P_{1}\xrightarrow{f_{1}}P_{0}\xrightarrow{f_{0}}\Gamma e_k\to \Gamma/I_k\to0
\end{equation}
with $P_{0},P_{1}\in\add\Gamma(1-e_k)$ satisfying the following conditions.
\begin{itemize}
\item[(a)] Applying $\Hom_\Gamma(-,\Gamma)$ to \eqref{projective resolution of Gamma/I}, we have an exact sequence
\begin{equation}\label{projective resolution of Gamma/I 2}
0\to e_k\Gamma\to\Hom_\Gamma(P_{0},\Gamma)\to\Hom_\Gamma(P_{1},\Gamma)\to e_k\Gamma\to \Gamma/I_k\to0.
\end{equation}
\item[(b)] Let $T_k^\ast=\Ker(T\otimes_\Gamma f_{0})$ and $\mu_k(T) =(T/T_k)\oplus T_k^{\ast}$.
Then $\mu_k(T)$ is a basic tilting $\Lambda$-module of projective dimension at most one and $T_k^{\ast} \not \simeq T_k$.
\item[(c)] If $Q_T$ and $Q_{\mu_k(T)}$ have no loops or 2-cycles, then $Q_{\mu_k(T)}\simeq\mu_k(Q_T)$ holds.
\end{itemize}
\end{proposition}

Notice that we can not drop the assumption $\Gamma/I_k\in\fl\Gamma$, which is automatically satisfied if $Q_T$ has no loops at $k$.
We need the following additional information for later application.

\begin{proposition}\label{3-CY mutation}
If $\Gamma/I_k\in\fl\Gamma$, then the following assertions hold.
\begin{itemize}
\item[(a)] There exist exact sequences (called \emph{exchange sequences})
\[0\to T_k^{\ast}\xrightarrow{g}U_k\xrightarrow{f}T_k\ \mbox{ and }\ 0\to T_k\xrightarrow{g'}U'_k\xrightarrow{f'}T^{\ast}_k\]
such that $f$ and $f'$ are right $\add (T/T_k)$-approximations, and $g$ and $g'$ are left $\add (T/T_k)$-approximations.
\item[(b)] In the category of finitely generated $\Lambda$-modules with projective dimension at most one,
$g$ and $g'$ are kernels of $f$ and $f'$ respectively, and $f$ and $f'$ are cokernels of $g$ and $g'$ respectively.
\item[(c)] The complex
$T_k\xrightarrow{g'}U'_k\xrightarrow{f'g}U_k\xrightarrow{f}T_k$
induces exact sequences
\begin{eqnarray*}
&0\to(T,T_k)\xrightarrow{g'}(T,U'_k)\xrightarrow{f'g}(T,U_k)\xrightarrow{f}[\add T/T_k](T,T_k)\to0,&\\
&0\to(T_k,T)\xrightarrow{f}(U_k,T)\xrightarrow{f'g}(U'_k,T)\xrightarrow{g'}[\add T/T_k](T_k,T)\to0.&
\end{eqnarray*}
If $Q_T$ has no loops at $k$, then we have $[\add T/T_k](T,T_k)=J_{\mod\Lambda}(T,T_k)$ and $[\add T/T_k](T_k,T)=J_{\mod\Lambda}(T_k,T)$.
\item[(d)] The complex
$T_k^{\ast}\xrightarrow{g}U_k\xrightarrow{fg'}U'_k\xrightarrow{f'}T_k^{\ast}$
induces exact sequences
\begin{eqnarray*}
&0\to(\mu_k(T),T_k^{\ast})\xrightarrow{g}(\mu_k(T),U_k)\xrightarrow{fg'}(\mu_k(T),U'_k)\xrightarrow{f'}[\add T/T_k](\mu_k(T),T_k^{\ast})\to0,&\\
&0\to(T_k^{\ast},\mu_k(T))\xrightarrow{f'}(U'_k,\mu_k(T))\xrightarrow{fg'}(U_k,\mu_k(T))\xrightarrow{g}[\add T/T_k](T_k^{\ast},\mu_k(T))\to0,&
\end{eqnarray*}
If $Q_T$ has no loops at $k$, then we have $[\add T/T_k](\mu_k(T),T_k^{\ast})=J_{\mod\Lambda}(\mu_k(T),T_k^{\ast})$
and $[\add T/T_k](T_k^{\ast},\mu_k(T))=J_{\mod\Lambda}(T_k^{\ast},\mu_k(T))$.
\end{itemize}
\end{proposition}

\begin{proof}
(a)(c) We have $\Tor^\Gamma_i(T,\Gamma/I_k)=0$ for any $i>1$.
Applying $T\otimes_\Gamma-$ to the sequence \eqref{projective resolution of Gamma/I}
and putting $U_k =T\otimes_\Gamma P_0$ and $U'_k = T \otimes_{\Gamma} P_1$,
we have the exact sequences in (a).

Since the functor $T\otimes_\Gamma\ \colon \add{}_\Gamma\Gamma\to\add{}_\Lambda T$ is an equivalence and
the sequence \eqref{projective resolution of Gamma/I} is exact,
we have that $f$ and $f'$ are right $\add (T/T_k)$-approximations, and that the upper sequence in (c) is exact.
Since the functor $\Hom_{\Gamma^{\op}}(-,T) \colon \add\Gamma_\Gamma\to\add{}_\Lambda T$ is a duality
and the sequence \eqref{projective resolution of Gamma/I 2} is exact,
we have that $g$ and $g'$ are left $\add (T/T_k)$-approximations, and that the lower sequence in (c) is exact.

(b) Clearly $g$ and $g'$ are kernels of $f$ and $f'$ respectively.
Since $\Lambda$ is 3-CY, we have $\Ext^i_\Lambda(\fl\Lambda,\Lambda)=0$ for $0\le i\le2$ \cite{ir}.
Thus any $\Gamma$-module $X$ with projective dimension at most one satisfies $\Ext^i_\Lambda(\fl\Lambda,X)=0$ for $i=0,1$.
On the other hand, since $\Gamma/I_k\in\fl\Gamma$ by our assumption,
we have that $\Cok f=T\otimes_\Gamma(\Gamma/I_k)$ and $\Cok f'=\Tor^\Gamma_1(T,\Gamma/I_k)$ are in $\fl\Lambda$.
Thus we have $\Ext^i_\Lambda(\Cok f,X)=0=\Ext^i_\Lambda(\Cok f',X)$ for $i=0,1$. This implies that $f$ and $f'$ are cokernels of $g$ and $g'$ respectively.

(d) Let $\Gamma'=\End_\Lambda(\mu_k(T))$ and $I'_k=\Gamma'(1-e_k)\Gamma'$.
We will show $\Gamma'/I'_k\in\fl\Gamma'$.
Then we will have the desired assertion by applying the argument in (c) for $\mu_k(T)$ instead of $T$.

For any $p\in J_{\mod\Lambda}(T_k^{\ast},T_k^{\ast})$,
there exists $q\in\End_\Lambda(U'_k)$ and $r\in\End_\Lambda(T_k)$ which make the diagram
\[\xymatrix{
0\ar[r]&T_k\ar^{g'}[r]\ar^r[d]&U'_k\ar^{f'}[r]\ar^q[d]&T_k^{\ast}\ar^p[d]\\
0\ar[r]&T_k\ar^{g'}[r]&U'_k\ar^{f'}[r]&T_k^{\ast}
}\]
commutative.
It is not difficult to check that the correspondence $p\mapsto r$
gives a well-defined isomorphism $\Gamma'/I'_k\to\Gamma/I_k$.

In particular, $Q_T$ has no loops at $k$ if and only if $Q_{\mu_k(T)}$ has no loops at $k$.
\end{proof}

\section{Cluster-tilted algebras are determined by their quivers}\label{sec2}

In this section we prove that cluster-tilted algebras are determined by their quivers.
Alternative proofs using potentials will be given as a consequence of the 
main result in section \ref{sec4}.

First we prove the following.

\begin{lemma}\label{lem2.1}
Let $T$ be a cluster-tilting object in some cluster category $\C_{Q}$, with $Q_T$ an acyclic
quiver. Then $\C_{Q_T}$ is equivalent to the cluster category $\C_Q$.
\end{lemma}

\begin{proof}
Although the assertion follows from the main result in \cite{kr2}, we give an elementary proof here.

We know that $T$ is induced by some tilting $KQ'$-module $U$, where $KQ'$ is derived equivalent to $KQ$.
Since $\End_{\C_{Q}}(T) \simeq KQ_T$, we have by Lemma \ref{l:acyc} that 
$\End_{KQ'}(U) \simeq \End_{\C_{Q'}}(T) \simeq \End_{\C_{Q}}(T) \simeq KQ_T$, 
so that $KQ'$ and $KQ_T$ are
also derived equivalent. Hence $KQ$ and $KQ_T$ are derived equivalent, so that the cluster categories $\C_Q$ and $\C_{Q_T}$ 
are equivalent.
\end{proof}

\begin{lemma}\label{lem2.2}
Let $Q$ be an acyclic quiver and $\C_Q$ the associated cluster category.
Let $T$ be a cluster-tilting object in $\C_Q$ such that the quiver $Q_T$ of $\End_{\C}(T)$ is isomorphic
to $Q$. Then there is an autoequivalence of triangulated categories
$F \colon \C_Q \to \C_Q$ with $F(T) \simeq KQ$.
\end{lemma}

\begin{proof}
We have that $T$ is induced by a tilting $KQ'$-module $U$, 
satisfying $\End_{KQ'}(U) \simeq KQ$ by Lemma \ref{l:acyc}.
Then we have a commutative diagram
$$
\xymatrix{
U \ar@{|->}[d] & \in & \D^{\rm b}(KQ') \ar[d] \ar^{\mathbf{R}\Hom_{KQ'}(U,-)}[rr] & & \D^{\rm b}(KQ) \ar[d] & \ni & KQ \ar@{|->}[d] \\
T              & \in  & \C_Q     \ar@{-->}^{F}[rr]          & & \C_Q           & \ni & KQ 
}
$$
where $\mathbf{R}\Hom_{KQ'}(U,-) \colon \D^{\rm b}(KQ') \to \D^{\rm b}(KQ)$ is an equivalence of triangulated categories
since $U$ is a tilting module \cite{h}. Hence $\mathbf{R}\Hom_{KQ'}(U,-)$ commutes with $\tau$ and $[1]$, so that
there is an equivalence $F \colon \C_Q \to \C_Q$ of triangulated categories,
with $F(T) \simeq KQ$. 
\end{proof}

Using this we get the following

\begin{theorem}\label{thm2.3}
Let $T_1$ and $T_2$ be cluster-tilting objects in the cluster categories $\C_{Q_1}$ and $\C_{Q_2}$ respectively, and assume
that the quivers $Q_{T_1}$ and $Q_{T_2}$ are isomorphic. Then 
there is an equivalence of triangulated categories $F \colon \C_{Q_1} \to \C_{Q_2}$ such that $F(T_1) \simeq T_2$.
In particular
$\End_{\C_{Q_1}}(T_1)\simeq\End_{\C_{Q_2}}(T_2)$ holds.
\end{theorem}

\begin{proof}
We have a sequence of mutations $\mu=\mu_{k_\ell}\circ\cdots\circ\mu_{k_1}$ such that $\mu(KQ_1)\simeq T_1$ in $\C_{Q_1}$.
Let $T$ be a cluster-tilting object in $\C_{Q_2}$ such that $\mu(T)\simeq T_2$.
Since $Q_{T_1} \simeq Q_{T_2}$, we have $Q_1 \simeq Q_T$. Hence by Lemmas \ref{lem2.1} and \ref{lem2.2}
there is an equivalence of triangulated categories $F \colon \C_{Q_1} \to \C_{Q_2}$ such that $F(KQ_1)\simeq T$. 
Since cluster-tilting mutations commute with equivalences of triangulated categories, we have
\[F(T_1)\simeq F(\mu(KQ_1)) \simeq \mu(F(KQ_1)) \simeq \mu(T)\simeq T_2,\]
and consequently $\End_{\C_{Q_1}} (T_1) \simeq \End_{\C_{Q_2}}(T_2)$.
\end{proof}

We can strengthen the above theorem as follows.

\begin{cor}
Let $T_1$ be a cluster-tilting object in a cluster category and $T_2$ a cluster-tilting object in an algebraic 
2-CY triangulated category $\C$ without loops or 2-cycles. If $Q_{T_1} \simeq Q_{T_2}$, then 
$\End_{\C}(T_1) \simeq \End_{\C}(T_2)$.
\end{cor}

\begin{proof}
By the same argument as in the proof of Theorem \ref{thm2.3}, there exists a sequence $\mu$ of mutations such that $\mu(T)\simeq T_2$ and $Q_T$ is acyclic.
Since $\C$ is assumed to be algebraic, it follows from \cite{kr2} that $\C$ is equivalent
to $\C_{Q_T}$, so the result follows from Theorem \ref{thm2.3}. 
\end{proof}

\section{Preliminaries on presentation of algebras}
In this section we consider presentations of algebras, which are important for the rest of this paper.

Let $K$ as before be an algebraically closed field and $\Gamma$ a $K$-algebra with
Jacobson radical $J_\Gamma$.
We regard $\Gamma$ as a topological algebra via the $J_\Gamma$-adic topology with a basic system $\{J_\Gamma^i\}_{i\ge0}$ 
of open neighborhoods of $0$.
Thus the closure of a subset $S$ of $\Gamma$ is given by
\begin{equation}\label{closure}
\overline{S}=\bigcap_{\ell\ge0}(S+J_\Gamma^\ell).
\end{equation}
We assume that $\Gamma$ satisfies the following conditions:
\begin{itemize}
\item[(A1)] $\Gamma/J_\Gamma\simeq\k^n$ for some $n >0$ and $\dim_{\k}(J_\Gamma/J_\Gamma^2)<\infty$,
\item[(A2)] the $J_\Gamma$-adic topology on $\Gamma$ is complete and separated, i.e., $\Gamma\simeq\varprojlim_{\ell\ge0}\Gamma/J_\Gamma^\ell$.
\end{itemize}
For example, finite dimensional $K$-algebras, complete path algebras $\widehat{KQ}$ of finite quivers $Q$,
and algebras over $R=K[[x_1,\ldots,x_d]]$ which are finitely generated $R$-modules satisfy (A1) and (A2).
If $\phi:\Lambda\to\Gamma$ is a homomorphism of $K$-algebras satisfying (A1), then $J_\Lambda\subset\phi^{-1}(J_\Gamma)$ holds.
In particular $\phi$ is continuous.

Let $Q$ be a finite quiver. It is convenient to introduce the following concepts.
For $a\in Q_1$, define a \emph{right derivative} $\partial^r_a \colon J_{\widehat{\k Q}}\to\widehat{\k Q}$
and a \emph{left derivative} $\partial^l_a \colon J_{\widehat{\k Q}}\to\widehat{\k Q}$ by
\begin{eqnarray*}
&\partial^r_a(a_1a_2\cdots a_{m-1}a_m)=
\left\{\begin{array}{cc}
a_1a_2\cdots a_{m-1}&\mbox{ if }\ a_m=a,\\
0&\mbox{ otherwise,}
\end{array}\right.&\\
&\partial^l_a(a_1a_2\cdots a_{m-1}a_m)=
\left\{\begin{array}{cc}
a_2\cdots a_{m-1}a_m&\mbox{ if }\ a_1=a,\\
0&\mbox{ otherwise,}
\end{array}\right.&
\end{eqnarray*}
and extend to $J_{\widehat{\k Q}}$ linearly and continuously. Then 
$\partial^r_a$ is a homomorphism of $\widehat{\k Q}$-modules and
$\partial^l_a$ is a homomorphism of $\widehat{\k Q}^{\op}$-modules.
Clearly for any $p\in J_{\widehat{\k Q}}$ we have
\[p=\sum_{a\in Q_1}(\partial^r_ap)a=\sum_{a\in Q_1}a(\partial^l_ap).\]
For simplicity we call $p\in\widehat{\k Q}$ \emph{basic} if $p$ is a formal linear sum of paths in $Q$ with a common 
start $i$ and a common end $j$.
In this case we write $s(p)=i$ and $e(p)=j$.
Clearly any ideal of $\widehat{\k Q}$ is generated by a set of basic elements.

We start with the following important fact.

\begin{proposition}\label{quiver generate}
Let $\Gamma$ be a $\k$-algebra satisfying (A1) and (A2), and let $Q$ be a finite quiver.
Assume that we have a $\k$-algebra homomorphism $\phi_0 \colon \k Q_0\to\Gamma$ and a homomorphism 
$\phi_1 \colon \k Q_1\to J_\Gamma$ of $(\k Q_0,\k Q_0)$-modules.
Then the following assertions hold.
\begin{itemize}
\item[(a)] $\phi_0$ and $\phi_1$ extend uniquely to a $\k$-algebra homomorphism 
$\phi \colon \widehat{\k Q}\to\Gamma$. Moreover $\Ker\phi$ is a closed ideal of $\widehat{\k Q}$.
\item[(b)] The following conditions are equivalent.
\begin{itemize}
\item[(i)] $\phi \colon \widehat{\k Q}\to\Gamma$ is surjective.
\item[(ii)] $\phi_0$ and $\phi_1$ induce surjections $\k Q_0\to \Gamma/J_\Gamma$ and $\k Q_1\to J_\Gamma/J_\Gamma^2$.
\item[(iii)] $\phi_0$ induces a surjection $\k Q_0\to \Gamma/J_\Gamma$, and for any $i\in Q_0$ we have an exact sequence
$${\displaystyle \bigoplus_{a\in Q_1,\ e(a)=i}\Gamma(\phi s(a))\xrightarrow{(\phi a)_a}J_\Gamma(\phi i)\to0.} $$
\item[(iv)]
$\phi_0$ induces a surjection $\k Q_0\to \Gamma/J_\Gamma$, and for any $i\in Q_0$ we have an exact sequence
$${\displaystyle \bigoplus_{a\in Q_1,\ s(a)=i}(\phi e(a))\Gamma\xrightarrow{(\phi a)_a}(\phi i)J_\Gamma\to0.}$$
\end{itemize}
\end{itemize}
\end{proposition}

\begin{proof}
(a) The first assertion is clear.
We observed that $\phi$ is continuous.
Since $\{0\}$ is a closed subset of $\Gamma$ then its inverse image $\Ker\phi$ is a closed subset of $\widehat{\k Q}$.

(b) (ii)$\Rightarrow$(i) $\phi$ induces a surjection $\phi_\ell\colon \k Q_\ell\to J_\Gamma^\ell/J_\Gamma^{\ell+1}$ for any $\ell$.
For any $x=x_0\in\Gamma$, there exist $p_\ell\in\k Q_\ell$ and $x_{\ell+1}\in J_\Gamma^{\ell+1}$
such that $x_{\ell+1}=x_\ell-\phi_\ell(p_\ell)$ for any $\ell\ge0$.
Then we have $\phi(\sum_{\ell\ge0}p_\ell)=x$.

\noindent (i)$\Rightarrow$(iii) 
It is enough to show that $\phi(p)$ is in the image of $(\phi a)_a$ for any
basic element $p\in J_{\widehat{KQ}}$ with $e(p)=i$.
This follows from the equality $p=\sum_{a\in Q_1}(\partial^r_ap)a$ since this implies $\phi p=(\phi\partial^r_ap)_a\cdot(\phi a)_a$.

\noindent (iii)$\Rightarrow$(ii) The sequence in (iii) induces a surjection
$\bigoplus_{a\in Q_1}\Gamma/J_\Gamma\xrightarrow{(\phi a)_a}J_\Gamma/J_\Gamma^2$.

\noindent The rest follows similarly.
\end{proof}

We shall need the following easy fact.

\begin{lemma}\label{closure of I}
Let $Q$ be a finite quiver and ${\bf S}$ a finite subset of $J_{\widehat{KQ}}$.
For the ideal $I=\langle{\bf S}\rangle$ of $\widehat{KQ}$, we have
$I=\sum_{v\in \mathbf{S}}(\widehat{\k Q})v+I\cdot Q_1$ and $\overline{I}=\sum_{v\in \mathbf{S}}(\widehat{\k Q})v+\overline{I}\cdot Q_1$.
\end{lemma}

\begin{proof}
Since $I$ is generated by $\mathbf{S}$, we have the first equality.

Let $p\in\overline{I}$.
For any $\ell > 0$ we can write $p=r_\ell+r'_\ell$ for some $r_\ell\in I$ and
$r'_\ell\in J_\Gamma^\ell$.
Then $p_\ell:=r_{\ell+1}-r_\ell=r'_\ell-r'_{\ell+1}\in I\cap
J_\Gamma^\ell$ for $\ell > 0$. Letting $p_0 =r_1$, we have
$p_0+p_1+\cdots+p_\ell=p-r'_{\ell+1}$. Thus $p=\sum_{\ell\ge0}p_\ell$.
We write $p_\ell=\sum_{v\in \mathbf{S}}p'_{\ell,v}v+\sum_{a\in Q_1}p''_{\ell,a}a$ for 
$p'_{\ell,v}\in\widehat{\k Q}$ and $p''_{\ell,a}\in I$.
Then $p=\sum_{v\in \mathbf{S}}(\sum_{\ell\ge0}p'_{\ell,v})v+\sum_{a\in Q_1}(\sum_{\ell\ge0}p''_{\ell,a})a$ holds.
Thus we have the second equality.
\end{proof}

We now give a result describing relations.

\begin{proposition}\label{quiver relation}
Assume that the conditions in Proposition \ref{quiver generate}(b) are satisfied.
For a finite set $\mathbf{S}$ of basic elements in $J_{\widehat{\k Q}}$\ , the following conditions are equivalent.
\begin{itemize}
\item[(a)] $\Ker\phi=\overline{I}$ holds for the ideal $I=\langle\mathbf{S}\rangle$ of $\widehat{\k Q}$.
\item[(b)] The following sequence is exact for any $i\in Q_0$.
\begin{eqnarray}\label{projective resolution}
\bigoplus_{v\in \mathbf{S},\ e(v)=i}\Gamma(\phi s(v))\xrightarrow{(\phi \partial^r_a v)_{v,a}}
\bigoplus_{a\in Q_1,\ e(a)=i}\Gamma(\phi s(a))\xrightarrow{(\phi a)_a}J_\Gamma(\phi i)\to0.
\end{eqnarray}
\item[(c)] The following sequence is exact for any $i\in Q_0$.
\begin{eqnarray*}
\bigoplus_{v\in \mathbf{S},\ s(v)=i}(\phi e(v))\Gamma\xrightarrow{(\phi \partial^l_a v)_{v,a}}
\bigoplus_{a\in Q_1,\ s(a)=i}(\phi e(a))\Gamma\xrightarrow{(\phi a)_a}(\phi i)J_\Gamma\to0.
\end{eqnarray*}
\end{itemize}
\end{proposition}

\begin{proof}
(a)$\Rightarrow$(b)
Since
$$(\phi\partial^r_av)_a\cdot(\phi a)_a=\phi(\sum_{a\in Q_1,\ e(a)=i}(\partial^r_av)a)=\phi v=0,$$
the sequence is a complex.
Now we assume that $(p_a)_a\in\bigoplus_{a\in Q_1\atop e(a)=i}\widehat{\k Q}(s(a))$ satisfies $(\phi p_a)_a\cdot(\phi a)_a=0$. Since
\[\sum_{a\in Q_1,\ e(a)=i}p_aa\in\Ker\phi=\overline{I},\]
there exists $q_v\in\widehat{\k Q}$ by Lemma \ref{closure of I} such that
\[\sum_{a\in Q_1,\ e(a)=i}p_aa-\sum_{v\in \mathbf{S}}q_vv\in\overline{I}\cdot Q_1.\]
Applying $\partial^r_a$ on both sides, we have
\[p_a-\sum_{v\in \mathbf{S}}q_v\partial^r_a v \in \overline{I}=\Ker\phi.\]
Thus $(\phi q_v)_v\in\bigoplus_{v\in \mathbf{S},\ e(v)=i}\Gamma(\phi s(v))$ satisfies
\[(\phi p_a)_a=(\phi q_v)_v\cdot(\phi\partial^r_av)_{v,a}.\]

(b)$\Rightarrow$(a) We shall show that $\Ker\phi=\overline{I}$.
Take any $v\in \mathbf{S}$ with $e(v)=i$. Since \eqref{projective resolution} is exact, we have
\[\phi v=\phi(\sum_{a\in Q_1,\ e(a)=i}(\partial^r_av)a)=(\phi\partial^r_av)_a\cdot(\phi a)_a=0.\]
This implies $v\in\Ker\phi$ and $I\subset\Ker\phi$.
Since $\Ker\phi$ is a closed ideal by Proposition \ref{quiver generate}, we have $\overline{I}\subset\Ker\phi$.

To prove $\Ker\phi\subset \overline{I}$, we will first show that, for any $p\in\Ker\phi$, 
there exists $p'\in I$ such that $p-p'\in (\Ker\phi)\cdot Q_1$.
Without loss of generality, we can assume that $p$ is basic with $e(p)=i$.
Since 
\[(\phi \partial^r_ap)_a\cdot(\phi a)_a=\phi p=0\]
holds, we have that $(\phi \partial^r_ap)_a$ factors through the left map $(\phi\partial^r_av)_{v,a}$ in \eqref{projective resolution}.
Thus there exists $q_v\in\widehat{\k Q}$ such that $(\phi \partial^r_ap)_a=(\phi q_v)_v\cdot(\phi \partial^r_av)_{v,a}$.
Then $$\partial^r_ap-\sum_{v\in \mathbf{S},\ e(v)=i}q_v \partial^r_av\in\Ker\phi.$$
Hence
\[p-\sum_{v\in \mathbf{S},\ e(v)=i}q_v v
=\sum_{a\in Q_1,\ e(a)=i}(\partial^r_ap-\sum_{v\in \mathbf{S},\ e(v)=i}q_v \partial^r_av)a\in(\Ker\phi)\cdot Q_1.\]
It follows that $p'=\sum_{v\in \mathbf{S},\ e(v)=i}q_v v\in I$ satisfies the desired condition.

Consequently we have $\Ker\phi\subset I+(\Ker\phi)\cdot Q_1$. This implies
\[\Ker\phi\subset I+(I+(\Ker\phi)\cdot Q_1)\cdot Q_1=I+(\Ker\phi)\cdot Q_1^2\subset\cdots\subset I+(\Ker\phi)\cdot Q_1^\ell\]
for any $\ell$. By \eqref{closure}, we have $\Ker\phi\subset\overline{I}$.

(a)$\Leftrightarrow$(c) can be shown dually.
\end{proof}

%

We give a relationship between the dimension of $\Ext^2$-spaces
and minimal sets of generators for relation ideals in complete path algebras (cf. \cite{bon,y}).
Note that the corresponding result is not true in general if we deal with ordinary path algebras $KQ$ instead of $\widehat{KQ}$.
See \cite[section 7]{bikr} for an example.
Note however that if $KQ/\langle{\bf S}\rangle$ is finite dimensional, then the corresponding result is true.

\begin{proposition}\label{p:rel}
Let $Q$ be a finite quiver, $J=J_{\widehat{KQ}}$\ ,
$I$ a closed ideal of $\widehat{KQ}$ contained in $J^2$, and $\Gamma = \widehat{KQ}/I$.
Let $\mathbf{S}$ be a finite set of basic elements in $I$.
\begin{itemize}
\item[(a)] $\mathbf{S}$ spans the $K$-vector space $I/(IJ+JI)$ if and only if $I=\overline{\langle\mathbf{S}\rangle}$.
Moreover $\mathbf{S}$ is a basis of the $K$-vector space $I/(IJ+JI)$ if and only if $\mathbf{S}$ is a minimal set satisfying $I=\overline{\langle\mathbf{S}\rangle}$.

\item[(b)] If $\mathbf{S}$ is a minimal set satisfying $I=\overline{\langle\mathbf{S}\rangle}$, then we have
$$\dim_K \Ext^2_{\Gamma}(S_j,S_i) = \#(i\widehat{KQ}j)\cap\mathbf{S} = \dim_K i(I/(IJ+JI))j$$
for any simple $\Gamma$-modules $S_i$ and $S_j$ associated with $i,j\in Q_0$.
\end{itemize}
\end{proposition}

\begin{proof}
(a) Assume first that $\mathbf{S}$ spans $I/(IJ+JI)$, and let $I'=\overline{\langle\mathbf{S}\rangle}$.
Then we have $I=I'+IJ+JI$. This implies
$$I=I'+(I'+IJ+JI)J+J(I'+IJ+JI)=I'+IJ^2+JIJ+J^2I=\cdots=I'+\sum_{0\le i\le\ell}J^iIJ^{\ell-i}$$
for any $\ell>0$. Thus we have $I\subset\bigcap_{\ell>0}(I'+J^\ell)=I'$ and hence $I=I'=\overline{\langle\mathbf{S}\rangle}$.

Next we assume that $I=\overline{\langle\mathbf{S}\rangle}$.
By Lemma \ref{closure of I}, we have $I=\langle\mathbf{S}\rangle+IJ$.
Thus $I=\langle\mathbf{S}\rangle+IJ+JI=K\mathbf{S}+IJ+JI$ holds.

The second equivalence is an immediate consequence of the first assertion.

(b) The minimality of $\mathbf{S}$ implies that the projective resolution in Proposition \ref{quiver relation}(b) is minimal. Thus we have the left equality.
The right equality follows from (a) since $(i\widehat{KQ}j)\cap\mathbf{S}$ is a basis of the $K$-vector space $i(I/(IJ+JI))j$.
\end{proof}

Now let $\C$ be a $\k$-linear category satisfying the following conditions:
\begin{itemize}
\item[(C1)] $\C$ is \emph{Krull-Schmidt}, i.e., any object in $\C$ is 
isomorphic to a finite direct sum of objects whose endomorphism rings are local.
\item[(C2)] $\End_{\C}(X)$ satisfies (A1) and (A2) for any basic object $X\in\C$.
\end{itemize}

For a finite quiver $Q$ and a field $K$, we denote by $\widehat{\k Q}$ also the category of 
finitely generated projective $\widehat{K Q}$-modules.
Then $\widehat{K Q}$ satisfies (C1) and (C2).
We have the following observation from Proposition \ref{quiver generate}(a).

\begin{lemma}\label{ext_cont}
For a category $\C$ satisfying (C1) and (C2) and a quiver $Q$, 
assume that an object $\Phi_0i\in \C$ for any $i\in Q_0$ and a morphism $\Phi_1a\in\C(\Phi s(a),\Phi e(a))$ for any $a\in Q_1$
are given. 
Then $\Phi_0$ and $\Phi_1$ extend uniquely to a $\k$-linear functor $\Phi \colon \widehat{\k Q}\to\C$.
\end{lemma}

Restating Proposition \ref{quiver relation}, we have the following result.

\begin{proposition}\label{category version}
Let $\C$ be a category satisfying (C1) and (C2), and let $Q$ be a finite quiver.
Let $\Phi \colon \widehat{\k Q}\to\C$ be a $\k$-linear functor, and let $T =\bigoplus_{i\in Q_0}\Phi i$.
For a finite set $\mathbf{S}$ of basic elements in $J_{\widehat{\k Q}}$\ , the following conditions are equivalent.
\begin{itemize}
\item[(a)] $\Phi$ induces a surjection $\phi \colon\widehat{\k Q}\to\End_{\C}(T)$ with $\Ker\phi=\overline{\langle \mathbf{S} \rangle}$.
\item[(b)] For any $i\in Q_0$, we have a complex
\[\bigoplus_{v\in \mathbf{S},\ e(v)=i}\Phi s(v)\xrightarrow{(\Phi \partial^r_av)_{v,a}}
\bigoplus_{a\in Q_1,\ e(a)=i}\Phi s(a)\xrightarrow{(\Phi a)_a}\Phi i\]
in $\add T$ which induces an exact sequence
\[\C(T,\bigoplus_{v\in \mathbf{S},\ e(v)=i}\Phi s(v))
\xrightarrow{(\Phi \partial^r_av)_{v,a}}\C(T,\bigoplus_{a\in Q_1,\ e(a)=i}\Phi s(a))
\xrightarrow{(\Phi a)_a}J_{\C}(T,\Phi i)\to0.\]
\item[(c)] For any $i\in Q_0$, we have a complex
\[\Phi i\xrightarrow{(\Phi a)_a}\bigoplus_{a\in Q_1,\ s(a)=i}\Phi e(a)
\xrightarrow{(\Phi \partial^l_av)_{a,v}}\bigoplus_{v\in \mathbf{S},\ s(v)=i}\Phi e(v)\]
in $\add T$ which induces an exact sequence
\[\C(\bigoplus_{v\in \mathbf{S},\ s(v)=i}\Phi e(v),T)
\xrightarrow{(\Phi \partial^l_av)_{a,v}}\C(\bigoplus_{a\in Q_1,\ s(a)=i}\Phi e(a),T)
\xrightarrow{(\Phi a)_a}J_{\C}(\Phi i,T)\to0.\]
\end{itemize}
\end{proposition}

\section{Jacobian algebras and weak 2-almost split sequences}\label{sec3}
In this section we study properties of Jacobian algebras of QP's.
We also investigate a variation of 2-almost split sequences and AR 4-angles
discussed in \cite{i,iy}, which we call weak 2-almost split sequences.
We show that there is a close relationship between Jacobian algebras and weak $2$-almost split sequences,
which will be used to prove the connection between mutation of QP's and of cluster-tilting objects in the next section. 

Let $(Q,W)$ be a QP.
For a pair of arrows $a$ and $b$, we define $\partial_{(a,b)}W$ by
\[\partial_{(a,b)}(a_1a_2\cdots a_m)=\sum_{a_i=a,\ a_{i+1}=b}a_{i+2}\cdots a_ma_1\cdots a_{i-1}\]
for any cycle $a_1\cdots a_m$ in $W$ and extend linearly and continuously.
Clearly $e(a)\neq s(b)$ implies $\partial_{(a,b)}W=0$.

The following easy observation is useful.

\begin{lemma}\label{partial*2}
\begin{itemize}
\item[(a)] $\sum_{a\in Q_1}(\partial_{(a,b)}W)a=\partial_bW=\sum_{c\in Q_1}c(\partial_{(b,c)}W)$ for any $b\in Q_1$.
\item[(b)] $\partial_{(a,b)}W=\partial^r_a(\partial_bW)=\partial^l_b(\partial_aW)$.
\end{itemize}
\end{lemma}

We have the following property of Jacobian algebras of QP's.

\begin{proposition}\label{projective complex}
Let $(Q,W)$ be a QP and $\Gamma=\P(Q,W)$.
Let $\phi:\widehat{KQ}\to\Gamma$ be a natural surjection.
Then there exist complexes
\begin{eqnarray*}
{\displaystyle\Gamma(\phi i)\xrightarrow{(\phi b)_b} \bigoplus_{\begin{smallmatrix}b\in Q_1\\ s(b)=i\end{smallmatrix}}\Gamma(\phi e(b))
\xrightarrow{(\phi \partial_{(a,b)}W)_{b,a}}\bigoplus_{\begin{smallmatrix}a\in Q_1\\ e(a)=i\end{smallmatrix}}\Gamma(\phi s(a))\xrightarrow{(\phi a)_a}J_\Gamma(\phi i)\to0,}\\
{\displaystyle(\phi i)\Gamma\xrightarrow{(\phi a)_a} \bigoplus_{\begin{smallmatrix}a\in Q_1\\ e(a)=i\end{smallmatrix}}(\phi s(a))\Gamma
\xrightarrow{(\phi \partial_{(a,b)}W)_{b,a}}\bigoplus_{\begin{smallmatrix}b\in Q_1\\ s(b)=i\end{smallmatrix}}(\phi e(b))\Gamma\xrightarrow{(\phi b)_b}(\phi i)J_\Gamma\to0}
\end{eqnarray*}
which are exact except at the second left terms.
\end{proposition}

\begin{proof}
Since $\partial_{(a,b)}W=\partial^r_a(\partial_bW)=\partial^l_b(\partial_aW)$ holds by Lemma \ref{partial*2},
the assertion follows immediately from Proposition \ref{quiver relation}(a)$\Rightarrow$(b)(c).
\end{proof}

As an immediate application, we have the following consequence.

\begin{proposition}\label{gorenstein}
Let $(Q,W)$ be a QP and $\Gamma =\P(Q,W)$.
\begin{itemize}
\item[(a)] $\Ext^2_\Gamma(S,\Gamma)=0=\Ext^2_{\Gamma^{\op}}(S',\Gamma)$ holds for any simple $\Gamma$-module $S$ and any simple $\Gamma^{\op}$-module $S'$.
\item[(b)] If $\Gamma$ is a finite dimensional $K$-algebra, then $\id{}_\Gamma\Gamma\le 1$ and $\id\Gamma_\Gamma\le 1$.
\end{itemize}
\end{proposition}

\begin{proof}
(b) follows immediately from (a).

(a) Let $S$ be a simple $\Gamma$-module. By Proposition \ref{projective complex}, there exists a complex
\[P\xrightarrow{f_2}P_1\xrightarrow{f_1}P_0\xrightarrow{f_0}P\]
of projective $\Gamma$-modules, which is exact at $P_0$, such that $\Cok f_0=S$ and
\begin{equation}\label{QQP}
\Hom_\Gamma(P_0,\Gamma)\xrightarrow{f_1}\Hom_\Gamma(P_1,\Gamma)\xrightarrow{f_2}\Hom_\Gamma(P,\Gamma)
\end{equation}
is exact.
We can take a morphism $f'_2 \colon P'\to P_1$ such that
\[P\oplus P'\xrightarrow{{f_2\choose f'_2}}P_1\xrightarrow{f_1}P_0\xrightarrow{f_0}P\to S\to0\]
gives a projective resolution of $S$.
Since \eqref{QQP} is exact, we have that
\[\Hom_\Gamma(P_0,\Gamma)\xrightarrow{f_1}\Hom_\Gamma(P_1,\Gamma)\xrightarrow{{f_2\choose f'_2}}\Hom_\Gamma(P\oplus P',\Gamma)\]
is exact. Consequently we have the desired properties.
\end{proof}


To study the property of Jacobian algebras of QP's categorically, we introduce the following concept.

\begin{definition}
Let $\C$ be a category satisfying (C1) and (C2), and $T\in\C$ an object.
We call a complex
\[U_1\xrightarrow{f_1}U_0\xrightarrow{f_0}X\]
in $\add T$ a \emph{right 2-almost split sequence} if
\begin{eqnarray*}
\C(T,U_1)\xrightarrow{f_1}\C(T,U_0)\xrightarrow{f_0}J_{\C}(T,X)\to0
\end{eqnarray*}
is exact.
In other words, $f_0$ is right almost split in $\add T$ and $f_1$ is a pseudo-kernel of $f_0$ in $\add T$.
Dually, we call a complex
\[X\xrightarrow{f_2}U_1\xrightarrow{f_1}U_0\]
in $\add T$ a \emph{left 2-almost split sequence} if
\begin{eqnarray*}
\C(U_0,T)\xrightarrow{f_1}\C(U_1,T)\xrightarrow{f_2}J_{\C}(X,T)\to0
\end{eqnarray*}
is exact. In other words, $f_2$ is left almost split in $\add T$ and $f_1$ is a pseudo-cokernel of $f_2$ in $\add T$.
We call a complex
\[X\xrightarrow{f_2}U_1\xrightarrow{f_1}U_0\xrightarrow{f_0}X\]
in $\add T$ a \emph{weak 2-almost split sequence} if
$U_1\xrightarrow{f_1}U_0\xrightarrow{f_0}X$ is a right 2-almost split sequence and
$X\xrightarrow{f_2}U_1\xrightarrow{f_1}U_0$ is a left 2-almost split sequence.
Note that we do not assume that $f_2$ (respectively, $f_0$) is a pseudo-kernel (respectively, pseudo-cokernel) of $f_1$.
\end{definition}

\noindent{\bf Remark } The definition of weak 2-almost split sequences
given here is more general than 2-almost split sequences and AR 4-angles in \cite{i,iy}.
One difference is that we assume neither right minimality of $f_0$ nor left minimality of $f_2$.
This point is necessary to deal with Jacobian algebras of non-reduced QP's.
Another difference is that our complex is not assumed to be a glueing of two exact sequences (respectively, triangles).
Note that we do not have uniqueness of weak 2-almost split sequences. 

\medskip
We give a key observation which gives a relationship between weak 2-almost split sequences and Jacobian algebras of QP's.

\begin{theorem}\label{key criterion}
Let $\C$ be a category satisfying (C1) and (C2), and let $(Q,W)$ be a QP. Let $\Phi \colon 
\widehat{\k Q}\to\C$ be a $\k$-linear functor,
and $T=\bigoplus_{i\in Q_0}\Phi i$. Then the following conditions are equivalent.
\begin{itemize}
\item[(a)] $\Phi$ induces an isomorphism $\P(Q,W)\simeq\End_{\C}(T)$.
\item[(b)] Any vertex $i$ in $Q$ gives rise to the following weak 2-almost split sequence in $\add T$.
\[{\displaystyle\Phi i\xrightarrow{(\Phi b)_b}\bigoplus_{b\in Q_1,\ s(b)=i}\Phi e(b)
\xrightarrow{(\Phi \partial_{(a,b)}W)_{b,a}}\bigoplus_{a\in Q_1,\ e(a)=i}\Phi s(a)\xrightarrow{(\Phi a)_a}\Phi i.}\]
\item[(c)] Any vertex $i$ in $Q$ gives rise to the following right 2-almost split sequence in $\add T$.
\[{\displaystyle\Phi e(b)\xrightarrow{(\Phi \partial_{(a,b)}W)_{b,a}}\bigoplus_{a\in Q_1,\ e(a)=i}\Phi s(a)\xrightarrow{(\Phi a)_a}\Phi i.}\]
\item[(d)] Any vertex $i$ in $Q$ gives rise to the following left 2-almost split sequence in $\add T$.
\[{\displaystyle\Phi i\xrightarrow{(\Phi b)_b}\bigoplus_{b\in Q_1,\ s(b)=i}\Phi e(b)
\xrightarrow{(\Phi \partial_{(a,b)}W)_{b,a}}\bigoplus_{a\in Q_1,\ e(a)=i}\Phi s(a).}\]
\end{itemize}
\end{theorem}

\begin{proof}
We shall apply Proposition \ref{category version} to $\mathbf{S} =\{\partial_aW\ |\ a\in Q_1\}$.

(b)$\Rightarrow$(c)(d) This is clear.

(c)$\Rightarrow$(a) By Lemma \ref{partial*2}(b), we have $\partial_{(a,b)}W=\partial^r_a(\partial_bW)$.
Looking at the complex
\[{\bigoplus_{b\in Q_1,\ s(b)=i}\Phi e(b)\xrightarrow{(\Phi \partial_{(a,b)}W)_{b,a}}
\bigoplus_{a\in Q_1,\ e(a)=i}\Phi s(a)\xrightarrow{(\Phi a)_a}\Phi i,}\]
we have the assertion by Proposition \ref{category version}(b)$\Rightarrow$(a).

(d)$\Rightarrow$(a) This is shown similarly.

(a)$\Rightarrow$(b) By Proposition \ref{category version}(a)$\Rightarrow$(b)(c), we have exact sequences
\begin{eqnarray*}
&{\displaystyle \C(T,\bigoplus_{b\in Q_1,\ s(b)=i}\Phi e(b))
\xrightarrow{(\Phi\partial^r_a(\partial_bW))_{b,a}}\C(T,\bigoplus_{a\in Q_1,\ e(a)=i}\Phi s(a))\xrightarrow{(a)_a}J_{\C}(T,\Phi i)\to0,}&\\
&{\displaystyle \C(\bigoplus_{a\in Q_1,\ e(a)=i}\Phi s(a),T)
\xrightarrow{(\Phi\partial^l_b(\partial_aW))_{b,a}}\C(\bigoplus_{b\in Q_1,\ s(b)=i}\Phi e(b),T)\xrightarrow{(b)_b}J_{\C}(\Phi i,T)\to0.}&
\end{eqnarray*}
Since $\partial_{(a,b)}W=\partial^r_a(\partial_bW)=\partial^l_b(\partial_aW)$ holds by Lemma \ref{partial*2}, we have the assertion.
\end{proof}

The following is an analogue of Theorem \ref{key criterion} which can be shown similarly and we will use later.

\begin{theorem}\label{analogue of key criterion}
Let $\C$ be a category satisfying (C1) and (C2), and let $(Q,W,F)$ be a QP with frozen vertices.
Let $\Phi \colon \widehat{\k Q}\to\C$ be a $\k$-linear functor,
and $T=\bigoplus_{i\in Q_0}\Phi i$. If $\Phi$ induces an isomorphism $\P(Q,W,F)\simeq\End_{\C}(T)$,
then any vertex $i$ in $Q_0\backslash F$ gives rise to the following weak 2-almost split sequence in $\add T$.
\[{\displaystyle\Phi i\xrightarrow{(\Phi b)_b}\bigoplus_{b\in Q_1,\ s(b)=i}\Phi e(b)
\xrightarrow{(\Phi \partial_{(a,b)}W)_{b,a}}\bigoplus_{a\in Q_1,\ e(a)=i}\Phi s(a)\xrightarrow{(\Phi a)_a}\Phi i.}\]
\end{theorem}


We now give the following sufficient condition for an
algebra to be a Jacobian algebra of a QP, which we will apply in section 6.
We call a cycle $a_1a_2\cdots a_m$ in a quiver $Q$ {\em full} if the vertices $s(a_1),\cdots,s(a_m)$ are distinct, and if there is an arrow $b$ in $Q$
with $s(b)$ and $e(b)$ in the cycle, then there is some arrow $a_i$ with $s(b)=s(a_i)$ and $e(b)=e(a_i)$.

\begin{proposition}\label{p:second}
Let $\Gamma$ be an algebra satisfying (A1) and (A2), and let $(Q,W)$
be a QP satisfying the following conditions.
\begin{itemize}
\item[(i)] There exists a surjective $K$-algebra homomorphism
$\phi \colon \P(Q,W)\to\Gamma$ such that $\Ker\phi$ is the closure of a finitely generated ideal.
\item[(ii)] Every cycle in $W$ is full.
\item[(iii)] The elements $\partial_aW$ for all arrows $a\in Q_1$ contained in cycles of $W$ are linearly independent over $K$.
\item[(iv)] $\dim_K\Ext^2_\Gamma(S,S')\le\dim_K\Ext^1_\Gamma(S',S)$
for any simple $\Gamma$-modules $S,S'$.
\end{itemize}
Then $\phi \colon \P(Q,W)\to\Gamma$ is an isomorphism.
\end{proposition}

\begin{proof}
Let $I$ be the kernel of the ring homomorphism
$\widehat{KQ}\to\Gamma$ induced by $\phi$ in (i). Then $\partial_aW\in
I$ for any $a\in Q_1$.

Let $J=J_{\widehat{KQ}}$.
We show that the image of $\partial_aW$ for all arrows $a\in Q_1$ contained in cycles of $W$ in $I/(IJ+JI)$ are linearly independent over $K$.
We assume $z=\sum_{i=1}^mc_i\partial_{a_i}W$ is in $IJ+JI$ for some
$c_1,\dots,c_m$ in $K$ and different arrows $a_1,\dots,a_m$ contained in cycles of $W$.
Fix a minimal set $\mathbf{S}'$ of basic elements satisfying $I=\overline{\langle\mathbf{S}'\rangle}$, which is finite by the assumption (i).
We have an equality $z=\sum_jp_jq_jr_j$ (possibly an infinite sum)
with $q_j\in\mathbf{S}'$ and $p_j,r_j\in\widehat{KQ}$ such that for
each $j$ one of $p_j$ and $r_j$ belongs to $J$.
Assume $z\neq0$. Then there exists a cycle $C$ in $W$ such
that $\partial_{a_i}C$ appears in $z$ for some $i=1,\dots,m$.
We have $\partial_{a_i}C=pqr$ for paths $p,q$ and $r$ appearing in
$p_j,q_j$ and $r_j$ respectively for some $j$, and one of $p$ and $r$
is a nontrivial path.
Let $S$ and $S'$ be simple $\Gamma$-modules corresponding to $s(q)$
and $e(q)$ respectively. Then we have $\Ext^2_\Gamma(S',S)\neq0$ by
applying Proposition \ref{p:rel}(b) to $\mathbf{S}'$.
By (iv), we have $\Ext^1_\Gamma(S,S')\neq0$, so there is an arrow in
$Q$ from $e(q)$ to $s(q)$.
Since $p$ or $r$ is a nontrivial path, it follows that the cycle $C$
is not full, a contradiction to (ii).
Consequently we have $z=0$. By (iii), we have that all $c_i$ are $0$.

For each arrow $j \overset{a}{\to} i$ lying on a cycle in $W$, we
have a relation of the form $\partial_a W$, hence $\dim_K \Ext^1_{\Gamma}(S_i,S_j)$
such relations, where $S_i$ and $S_j$ denote the simple modules at the vertices $i$ and $j$. 
We then have
\begin{multline*}
\dim_K \Ext^1_{\Gamma}(S_i,S_j) \leq \dim_K i(I/(IJ+JI))j \\ 
= \dim_K \Ext^2_{\Gamma}(S_j,S_i) \leq \dim_K \Ext^1_{\Gamma}(S_i,S_j),
\end{multline*}
where the first inequality follows by the above, the equality follows by Proposition \ref{p:rel}(b), and the 
last inequality holds by (iv). Hence we have $I=\overline{\langle\partial_aW\ |\ a\in Q_1\rangle}$ by Proposition \ref{p:rel}(a).
It follows that $\phi$ is an isomorphism.
\end{proof} 

Now we study the relationship between weak 2-almost split sequences and exchange triangles/sequences for the following three cases.
\begin{itemize}
\item[(A)] $\C=\mod\Lambda$ for a 3-CY algebra $\Lambda$ and $T$ is a basic tilting $\Lambda$-module of projective dimension at most one
(section \ref{subsec:mut_obj2}).
\item[(B)] $\C$ is a 2-CY Frobenius category and $T$ is a basic cluster-tilting object.
\item[(C)] $\C$ is a 2-CY triangulated category and $T$ is a basic cluster-tilting object.
\end{itemize}
We also consider the following stronger version of (C) associated with $k\in Q_0$.
\begin{itemize}
\item[(C$_k$)] In (C), we have $\Hom_\Lambda(\Ext^1_\Lambda(D\Lambda,S_k),S_k)=0$ for $\Lambda:=\End_{\C}(T)$
and the simple $\Lambda$-module $S_k$. (We call this the \emph{vanishing condition} at $k$.)
\end{itemize}

The following observation shows that weak $2$-almost split sequences appear naturally in cluster-tilting theory.

\begin{proposition}\label{gluing of exchange}
Let $\C$ and $T=T_1\oplus\cdots\oplus T_n$ satisfy one of the above (A), (B) or (C).
Assume that the quiver of $\End_{\C}(T)$ has no loops at $k$.
Glueing exchange triangles/sequences
\[T_k^{\ast}\xrightarrow{g}U_k\xrightarrow{f}T_k\ \mbox{ and }\ T_k\xrightarrow{g'}U'_k\xrightarrow{f'}T_k^{\ast}\]
given in sections \ref{subsec:mut_obj} and \ref{subsec:mut_obj2},
we have weak 2-almost split sequences
\begin{eqnarray*}
&T_k\xrightarrow{g'}U'_k\xrightarrow{f'g}U_k\xrightarrow{f}T_k\ \mbox{ in }\ \add T,&\\
&T_k^{\ast}\xrightarrow{g}U_k\xrightarrow{fg'}U'_k\xrightarrow{f'}T_k^{\ast}\ \mbox{ in }\ \add\mu_k(T).&
\end{eqnarray*}
\end{proposition}

\begin{proof}
The cases (B) and (C) are easily shown (cf. \cite[Theorem 3.10]{iy}), and the case (A) is shown in Proposition \ref{3-CY mutation}(c)(d).
\end{proof}

In the rest of this section, we shall show the following converse statement
of Proposition \ref{gluing of exchange}, which plays an important role in the next section.

\begin{theorem}\label{from 2-slmost split to exchange}
Let $\C$ and $T=T_1\oplus\cdots\oplus T_n$ satisfy one of the above (A), (B) or (C$_k$).
Assume that $Q_T$ has no loops at $k$. Then for any weak $2$-almost split sequence
\begin{eqnarray*}
&T_k\xrightarrow{f_2}U_1\xrightarrow{f_1}U_0\xrightarrow{f_0}T_k&
\end{eqnarray*}
with $f_1\in J_{\C}$, there exist exchange triangles/sequences
\begin{eqnarray*}
&T_k^{\ast}\xrightarrow{g}U_0\xrightarrow{f_0}T_k\ \ \mbox{ and }\ \ 
T_k\xrightarrow{f_2}U_1\xrightarrow{g'}T_k^{\ast}&
\end{eqnarray*}
such that $f_1=g'g$.
\end{theorem}

We first give a proof for the case (A).
Since $f_1\in J_{\C}$, we have that $f_0$ and $f_2$ are minimal right and left $\add(T/T_k)$-approximations respectively.
Consider \emph{any} exchange sequences
\begin{eqnarray*}
&0\to T_k^{\ast}\xrightarrow{g}U_0\xrightarrow{f_0}T_k\ \ \mbox{ and }\ \ 
0\to T_k\xrightarrow{f_2}U_1\xrightarrow{g'}T_k^{\ast}.&
\end{eqnarray*}
There is an automorphism $p\in{\rm Aut}_{\C}(U_1)$ such that $f_1=pg'g$.
Since $(f_2pg')g=f_2f_1=0$ and $g$ is injective, we have $f_2pg'=0$.
Thus there exists $q\in{\rm Aut}_{\C}(T_k)$ such that $qf_2=f_2p$.
\[\xymatrix{
&T_k\ar^{f_2}[r]\ar^q[d]&U_1\ar^{f_1}[r]\ar^p[d]&U_0\ar^{f_0}[r]\ar@{=}[d]&T_k\ar@{=}[d]\\
0\ar[r]&T_k\ar^{f_2}[r]&U_1\ar^{g'g}[r]&U_0\ar^{f_0}[r]&T_k
}\]
Since $f_2=qf_2p^{-1}$, we have exchange sequences
\begin{eqnarray*}
&0\to T_k^{\ast}\xrightarrow{g}U_0\xrightarrow{f_0}T_k\ \ \mbox{ and }\ \ 
0\to T_k\xrightarrow{qf_2p^{-1}=f_2}U_1\xrightarrow{pg'}T_k^{\ast}&
\end{eqnarray*}
satisfying $pg'g=f_1$.

\medskip 
The proof for the case (B) is the same.

\medskip 
In the rest of this section, we consider the case (C$_k$).
We need preliminary results.

\begin{lemma}\label{Ext^1 top}
Let $\Lambda$ be a basic finite dimensional $K$-algebra with $\pd{}_\Lambda(D\Lambda)\le 1$.
If the quiver of $\Lambda$ has no loops at $k$ and $\Lambda$ satisfies the vanishing condition at $k$, then $P=\Lambda e_k$ satisfies
\[\Ext^1_\Lambda(\nu P,P)=
\Hom_\Lambda(\nu P,\nu(\Lambda/P))\Ext^1_\Lambda(\nu(\Lambda/P),P)+
\Ext^1_\Lambda(\nu P,\Lambda/P)\Hom_\Lambda(\Lambda/P,P).\]
\end{lemma}

\begin{proof}
Let $M=\Ext^1_\Lambda(D\Lambda,\Lambda)$ and $e=e_k$.
We have $\Ext^1_\Lambda(D\Lambda,S_k)=M\otimes_\Lambda S_k$ since $\pd{}_\Lambda(D\Lambda)\le 1$.
The vanishing condition at $k$ implies
\[DS_k\otimes_\Lambda M\otimes_\Lambda S_k=D\Hom_\Lambda(M\otimes_\Lambda S_k,S_k)=0,\]
so we have $e(\top_{\Lambda\otimes_{\k}\Lambda^{\rm op}}M)e=0$. This implies
\[eMe=e(J_\Lambda M+MJ_\Lambda)e.\]
Since the quiver of $\Lambda$ has no loops at $k$, we have $eJ_\Lambda=e\Lambda(1-e)\Lambda$ and $J_\Lambda e=\Lambda(1-e)\Lambda e$.
Thus we have
\[eMe=e(J_\Lambda M+MJ_\Lambda)e=e\Lambda(1-e)Me+eM(1-e)\Lambda e.\]
Since $(D\Lambda)(1-e)=\nu(\Lambda/P)$, we have
\[(1-e)Me=\Ext^1_\Lambda((D\Lambda)(1-e),\Lambda e)=\Ext^1_\Lambda(\nu(\Lambda/P),P).\]
Similarly, we have $eM(1-e)=\Ext^1_\Lambda(\nu P,\Lambda/P)$.
Thus we have the desired equality.
\end{proof}

We also need the following easy observation.

\begin{lemma}\label{C and Ext^1}
Let $\C$ be a 2-CY triangulated category with a cluster-tilting object $T$.
Let $\Lambda=\End_{\C}(T)$ and $F=\C(T,-):\C\to\mod\Lambda$.
Then we have a functorial isomorphism
\[\alpha_{X,Y}:\C(X[1],Y)\simeq\Ext^1_\Lambda(\nu F(X),F(Y))\]
for any $X,Y\in\add T$.
\end{lemma}

\begin{proof}
We have an equivalence $F:\add T\xrightarrow{\sim}\add\Lambda$.
Take a triangle
\begin{equation}\label{approximation}
X[1]\to U_1\to U_0\to X[2]
\end{equation}
with $U_0,U_1\in\add T$.
Applying $F$ and using $F(X[2])\simeq D\C(X,T)=\nu F(X)$, we have a projective resolution
\begin{equation}\label{(T,approximation)}
0\to F(U_1)\to F(U_0)\to\nu F(X)\to0
\end{equation}
of the injective $\Lambda$-module $\nu F(X)$.

Applying $\C(-,Y)$ to \eqref{approximation} and
$\Hom_\Lambda(-,F(Y))$ to \eqref{(T,approximation)}
and comparing them by Yoneda's Lemma on $\add T$,
we have a commutative diagram
\[\xymatrix@C=1em{
\C(U_0,Y)\ar[r]\ar^{\wr}[d]&\C(U_1,Y)\ar[r]\ar^{\wr}[d]&\C(X[1],Y)\ar[d]\ar[r]&
\C(U_0[-1],Y)=0\\
\Hom_\Lambda(F(U_0),F(Y))\ar[r]&\Hom_\Lambda(F(U_1),F(Y))\ar[r]&
\Ext^1_\Lambda(\nu F(X),F(Y))\ar[r]&0}\]
of exact sequences. Thus we have an isomorphism
\begin{equation*}
\C(X[1],Y)\simeq\Ext^1_\Lambda(\nu F(X),F(Y)),
\end{equation*}
which is easily checked to be functorial for $X,Y\in\add T$.
\end{proof}

Now we have the following result.

\begin{proposition}\label{T[1] to T}
In the case (C$_k$), any morphism in $\C(T_k[1],T_k)$ factors through $\add((T/T_k)[1]\oplus(T/T_k))$.
\end{proposition}

\begin{proof}
We use the notation in Lemma \ref{C and Ext^1}.
Let $P=F(T_k)$ and $I=\nu P$. Then we have an isomorphism
\[\alpha_{T_k,T_k}:\C(T_k[1],T_k)\simeq\Ext^1_\Lambda(I,P).\]
It is easily checked that we have a commutative diagram
\[\xymatrix@C=1.5cm{
\C(T_k[1],T/T_k)\times\C(T/T_k,T_k)\ar^{\alpha_{T_k,T/T_k}\times F}_{\sim}[r]\ar^{\rm comp.}[d]&
\Ext^1_\Lambda(I,\Lambda/P)\times\Hom_\Lambda(\Lambda/P,P)\ar[d]\\
\C(T_k[1],T_k)\ar^{\alpha_{T_k,T_k}}_{\sim}[r]&\Ext^1_\Lambda(I,P)}\]
whose horizontal maps are isomorphisms.
Comparing the images of vertical maps, we have that $\alpha_{T_k,T_k}$ induces an isomorphism
\begin{equation}\label{T/T_k}
[T/T_k](T_k[1],T_k)\simeq\Ext^1_\Lambda(I,\Lambda/P)\Hom_\Lambda(\Lambda/P,P).
\end{equation}
Similarly, we have a commutative diagram
\[\xymatrix@C=2cm{
\C(T_k[1],(T/T_k)[1])\times\C((T/T_k)[1],T_k)\ar^{F[1]\times\alpha_{T/T_k,T_k}}_{\sim}[r]\ar^{\rm comp.}[d]&
\Hom_\Lambda(I,\nu(\Lambda/P))\times\Ext^1_\Lambda(\nu(\Lambda/P),P)\ar[d]\\
\C(T_k[1],T_k)\ar^{\alpha_{T_k,T_k}}_{\sim}[r]&\Ext^1_\Lambda(I,P)}\]
whose horizontal maps are isomorphisms.
Comparing the images of vertical maps, we have that $\alpha_{T_k,T_k}$ induces an isomorphism
\begin{equation}\label{T/T_k[1]}
[(T/T_k)[1]](T_k[1],T_k)\simeq
\Hom_\Lambda(I,\nu(\Lambda/P))\Ext^1_\Lambda(\nu(\Lambda/P),P).
\end{equation}
Since $Q_T$ has no loops at $k$, we have
\[\Ext^1_\Lambda(I,P)=\Hom_\Lambda(I,\nu(\Lambda/P))\Ext^1_\Lambda(\nu(\Lambda/P),P)
+\Ext^1_\Lambda(I,\Lambda/P)\Hom_\Lambda(\Lambda/P,P)\]
by Lemma \ref{Ext^1 top}. Using \eqref{T/T_k} and \eqref{T/T_k[1]}, we have
\[\C(T_k[1],T_k)=[(T/T_k)[1]](T_k[1],T_k)+[T/T_k](T_k[1],T_k),\]
which shows the assertion.
\end{proof}

Now we are ready to prove Theorem \ref{from 2-slmost split to exchange} for the case (C$_k$).

Since $f_1\in J_{\C}$, we have that $f_0$ and $f_2$ are minimal right and left $\add(T/T_k)$-approximations respectively.
Consider \emph{any} exchange triangles 
\begin{eqnarray*}
&T_k^{\ast}\xrightarrow{g}U_0\xrightarrow{f_0}T_k\xrightarrow{e}T_k^{\ast}[1]\ \ \mbox{ and }\ \ 
T_k\xrightarrow{f_2}U_1\xrightarrow{g'}T_k^{\ast}\xrightarrow{e'}T_k[1].&
\end{eqnarray*}
There are automorphisms $p\in{\rm Aut}_{\C}(U_1)$ and $q\in{\rm Aut}_{\C}(U_0)$
such that $f_1=pg'g=g'gq$.
If $g'g=0$, then $f_1=0$ and the assertion holds. Thus we can assume $g'g\neq0$.

(i) First we will show that there exists $u\in\End_{\C}(T_k^{\ast})$ such that $g'ug=f_1$.

Consider the following commutative diagram of triangles
\[\xymatrix{
U_1\ar^{g'}[r]\ar_{pg'}[d]&T_k^{\ast}\ar^{e'}[r]\ar_{gq}[d]&T_k[1]\ar^{f_2[1]}[r]&
U_1[1]\ar^{(pg')[1]}[d]\\
T_k^{\ast}\ar^{g}[r]&U_0\ar^{f_0}[r]&T_k\ar^e[r]&T_k^{\ast}[1]
}\]
We have a morphism $h\in\C(T_k[1],T_k)$ which keeps the above diagram commutative.
By Proposition \ref{T[1] to T}, there exist $s\in\C(U_1[1],T_k)$ and $t\in\C(T_k[1],U_0)$ such that
$h=f_2[1]s+tf_0$.
Since $(gq-e't)f_0=e'(h-tf_0)=e'(f_2[1])s=0$,
there exists $u\in\C(T_k^{\ast},T_k^{\ast})$ such that $gq=e't+ug$.
\[\xymatrix{
U_1\ar^{g'}[r]\ar_{pg'}[d]&T_k^{\ast}\ar^{e'}[r]\ar_{gq}[d]\ar^{u}[dl]&
T_k[1]\ar^{f_2[1]}[r]\ar^{t}[dl]\ar^h[d]&U_1[1]\ar^{(pg')[1]}[d]\ar^{s}[dl]\\
T_k^{\ast}\ar^{g}[r]&U_0\ar^{f_0}[r]&T_k\ar^e[r]&T_k^{\ast}[1]
}\]
In particular, we have $f_1=g'gq=g'ug$.

(ii) Next we will show that $u$ is an automorphism.
Then we have the assertion since we have exchange triangles
\begin{eqnarray*}
&T_k^{\ast}\xrightarrow{g}U_0\xrightarrow{f_0}T_k\xrightarrow{e}T_k^{\ast}[1]\ \ \mbox{ and }\ \ 
T_k\xrightarrow{f_2}U_1\xrightarrow{g'u}T_k^{\ast}\xrightarrow{u^{-1}e'}T_k[1]&
\end{eqnarray*}
satisfying $f_1=(g'u)g$.

Assume $u\in J_{\C}$. Then there exists
$v\in J_{\C}(U_1,U_1)$ such that $g'u=vg'$.
Since $pg'g=g'ug=vg'g$, we have $(p-v)g'g=0$.
Since $p-v$ is an automorphism by $p\in{\rm Aut}_{\C}(U_1)$ and $v\in J_{\C}(U_1,U_1)$, we have $g'g=0$.
This contradicts our assumption.
\qed

\section{Mutation of quivers with potentials and cluster-tilting objects}\label{sec4}

\subsection{Main results}

In this section we use the results from section \ref{sec3} to show that there is a close connection between
mutation of quivers with potentials on one hand and mutation of cluster-tilting objects in 2-CY triangulated categories in section \ref{subsec:mut_obj},
or mutation of tilting modules over 3-CY algebras in section \ref{subsec:mut_obj2}, on the other hand. More specifically 
we will prove the following main results in this section.

The first main result is about mutation of tilting modules over 3-CY algebras.

\begin{theorem}\label{main2}
Let $\Lambda$ be a 3-CY algebra with a basic tilting module $T$ of projective dimension at most one.
If $\End_\Lambda(T)\simeq\P(Q,W)$ for a QP $(Q,W)$ and no 2-cycles start in the vertex $k$ of $Q$, then $\End_\Lambda(\mu_k(T))\simeq\P(\mu_k(Q,W))$.
\end{theorem}

The second and third main results are about mutation of cluster-tilting objects in 2-CY triangulated categories.

Recall that a finite dimensional algebra $\Lambda$ satisfies \emph{the vanishing condition at $k$}
if $\Hom_\Lambda(\Ext^1_\Lambda(D\Lambda,S_k),S_k)=0$ holds for the simple $\Lambda$-module $S_k$.

\begin{theorem}\label{main}
Let $\C$ be a 2-CY triangulated category with a basic cluster-tilting object $T$.
If $\End_{\C}(T)\simeq \P(Q,W)$ for a QP $(Q,W)$, no 2-cycles start in the vertex $k$ of $Q$
and the vanishing condition (or more generally, the glueing condition below) is satisfied at $k$,
then $\End_{\C}(\mu_k(T))\simeq\P(\mu_k(Q,W))$.
\end{theorem}

Here the \emph{glueing condition at $k$} means that there exist exchange triangles for $T_k$ such that
their glueing is isomorphic to the weak 2-almost split sequence
\[{\displaystyle\Phi k\xrightarrow{(\Phi b)_b}\bigoplus_{b\in Q_1,\ s(b)=k}\Phi e(b)
\xrightarrow{(\Phi \partial_{(a,b)}W)_{b,a}}\bigoplus_{a\in Q_1,\ e(a)=k}\Phi s(a)\xrightarrow{(\Phi a)_a}\Phi k}\]
given in Theorem \ref{key criterion}(b).
By Theorem \ref{from 2-slmost split to exchange}, we know that the vanishing condition at $k$ implies the glueing condition at $k$.

\medskip
To state the last main result, we prepare some terminology.

Let $\C$ be a 2-CY triangulated category with a basic cluster-tilting object $T$.
We say that a homomorphism $\phi:\P(Q,W)\to\End_{\C}(T)$ of $\k$-algebras is \emph{liftable} if the following conditions are satisfied.
\begin{itemize}
\item There exists a Frobenius category $\E$ satisfying (C1) and (C2) such that $\C$ is triangle equivalent to the stable category of $\E$.
\item There exists QP $(\widetilde{Q},\widetilde{W},F)$ with frozen vertices such that
$(Q,W)$ is the QP obtained from $(\widetilde{Q},\widetilde{W},F)$ by removing all vertices in $F$.
\item There exists a $\k$-linear functor $\Phi:\widehat{\k\widetilde{Q}}\to\E$ such that
$\bigoplus_{i\in Q_0}\Phi i$ is $T$ in $\C$, $\bigoplus_{i\in F}\Phi i$ is a progenerator of $\E$
and $\phi$ is induced from $\Phi$.
\item For any $i\in Q_0$, we have the following weak 2-almost split sequence in $\add T\subset\E$.
\[{\displaystyle\Phi i\xrightarrow{(\Phi b)_b}\bigoplus_{b\in \widetilde{Q}_1,\ s(b)=i}\Phi e(b)
\xrightarrow{(\Phi \partial_{(a,b)}\widetilde{W})_{b,a}}\bigoplus_{a\in \widetilde{Q}_1,\ e(a)=i}\Phi s(a)\xrightarrow{(\Phi a)_a}\Phi i.}\]
\end{itemize}
In this case, it is easily checked that $\phi:\P(Q,W)\to\End_{\C}(T)$ is an isomorphism.

\begin{theorem}\label{main3}
Let $\C$ be a 2-CY triangulated category with a basic cluster-tilting object $T$.
If we have a liftable isomorphism $\End_{\C}(T)\simeq \P(Q,W)$ for a QP $(Q,W)$ and no 2-cycles start in the vertex $k$ of $Q$,
then we have a liftable isomorphism $\End_{\C}(\mu_k(T))\simeq\P(\mu_k(Q,W))$.
\end{theorem}

Now we pose the following.

\medskip\noindent
{\bf Conjecture }
\emph{Theorem \ref{main} is true without assuming the vanishing condition.}

\medskip
We have the following direct consequences of Theorems \ref{main2}, \ref{main} and \ref{main3}.

\begin{cor}\label{c:mut-class}
Let $\C$ be a $\k$-linear category and $T\in\C$ be an object such that
we have an isomorphism $\phi:\End_{\C}(T)\simeq \P(Q,W)$ of $\k$-algebras for a QP $(Q,W)$.
Let $k_1,\dots,k_\ell$ be a sequence of vertices in $Q$ such that
$k_i$ does not lie on a 2-cycle in $\mu_{k_{i-1}}\circ\cdots\circ\mu_{k_1}(Q,W)$ for each $i=1,\ldots,\ell$.
If one of the following conditions is satisfied, then we have
$$\End_{\C}(\mu_{k_\ell}\circ\cdots\circ\mu_{k_1}(T))\simeq \P(\mu_{k_\ell}\circ\cdots\circ\mu_{k_1}(Q,W)).$$
\begin{itemize}
\item[(a)] $\C=\mod\Lambda$ for a 3-CY algebra $\Lambda$ and $T$ is a basic tilting $\Lambda$-module of projective dimension at most one.
\item[(b)] $\C$ is a 2-CY triangulated category, $T$ is a basic cluster-tilting object and $\phi$ is liftable.
\item[(c)] $\C$ is a 2-CY triangulated category, $T$ is a basic cluster-tilting object
and $\End_{\C}(\mu_{k_{i-1}}\circ\cdots\circ\mu_{k_1}(T))$ satisfies the vanishing condition at $k_i$ for each $i=1,\ldots,\ell$.
\end{itemize}
\end{cor}

It is not known in general whether the algebra $\End_{\C}(T^{\ast})$
can be defined directly from the algebra $\End_{\C}(T)$ when
$T^{\ast}$ is a cluster-tilting object in the triangulated 2-CY category $\C$
obtained from the cluster tilting object $T$ by mutation.
When the potential $(Q',W')$ is obtained by mutation from the potential $(Q,W)$,
it is also not known if the Jacobian algebra $\P(Q',W')$ can be defined directly
from the Jacobian algebra $\P(Q,W)$.
However, we have the following.

\begin{theorem}\label{unique}
Let $\Lambda$ be an algebra isomorphic to a Jacobian algebra $\P(Q,W)$ for some QP $(Q,W)$ and
to a 2-CY-tilted algebra $\End_{\C}(T)$ for some cluster-tilting object $T$ in a 2-CY category $\C$.
Assume $k$ is a vertex in $Q$ not lying on any 2-cycle.
Assume that the isomorphism $\P(Q,W)\to\End_{\C}(T)$ is liftable,
or $\Lambda$ satisfies the vanishing condition at $k$.
\begin{itemize}
\item[(a)] The Jacobian algebra $\P(\mu_k(Q,W))$ is determined by the algebra $\Lambda$ and does not depend on the choice of a QP $(Q,W)$.
\item[(b)] The 2-CY-tilted algebra $\End_{\C}(\mu_k(T))$ is determined by the algebra $\Lambda$ and does not depend on the choice of 
a 2-CY category $\C$ and a cluster-tilting object $T$.
\end{itemize}
\end{theorem}

\begin{proof}
This is a direct consequence of Theorems \ref{main} and \ref{main3}.
\end{proof}

In section 6, we discuss an interesting class of triangulated 2-CY categories where we have liftable isomorphisms
for a large class of 2-CY-tilted algebras.

\subsection{Proof of main results}

To prove Theorems \ref{main2}, \ref{main} and \ref{main3} at the same time, we start with the following general setup:

Let $\C$ be a $\k$-linear category satisfying (C1) and (C2) and $(Q,W,F)$ be a reduced QP with frozen vertices.
Assume that we have a $\k$-linear functor $\Phi:\widehat{\k Q}\to\C$. Let $T_i=\Phi i$ and $T=\bigoplus_{i\in Q_0}T_i$,
and we simply denote $\Phi p$ by $p$ for any morphism $p$ in $\widehat{KQ}$.
We assume the following condition is satisfied.
\begin{itemize}
\item[(O)] For any $i\in Q_0\backslash F$, we have a weak 2-almost split sequence
\[{\displaystyle T_i\xrightarrow{(b)_b}\bigoplus_{b\in Q_1,\ s(b)=i}T_{e(b)}
\xrightarrow{(\partial_{(a,b)}W)_{b,a}}\bigoplus_{a\in Q_1,\ e(a)=i}T_{s(a)}\xrightarrow{(a)_a}T_i}\]
in $\add T$, which we simply denote by
\begin{eqnarray*}
T_i\xrightarrow{f_{i2}}U_{i1}\xrightarrow{f_{i1}}U_{i0}\xrightarrow{f_{i0}}T_i.
\end{eqnarray*}
\end{itemize}
Moreover we fix $k\in Q_0\backslash F$
and assume that there exists an indecomposable object $T_k^{\ast}\in\C$ which
does not belong to $\add T$ such that the following conditions are satisfied, where we let $\mu_k(T) =(T/T_k)\oplus T_k^{\ast}$.
\begin{itemize}
\item[(I)] There exist complexes
\begin{eqnarray*}
&T_k\xrightarrow{f_{k2}}U_{k1}\xrightarrow{h_k}T_k^{\ast}\ \ \mbox{ and }\ \ T_k^{\ast}\xrightarrow{g_k}U_{k0}\xrightarrow{f_{k0}}T_k&
\end{eqnarray*}
in $\C$ such that $f_{k1}=h_kg_k$.
\item[(II)] We have the following weak 2-almost split sequence in $\add\mu_k(T)$.
\[T_k^{\ast}\xrightarrow{g_k}U_{k0}\xrightarrow{f_{k0}f_{k2}}U_{k1}\xrightarrow{h_k}T_k^{\ast}\]
\item[(III)] The following sequences are exact.
\begin{eqnarray*}
&(T_k^{\ast},T_k^{\ast})\xrightarrow{g_k}(T_k^{\ast},U_{k0})\xrightarrow{f_{k0}}(T_k^{\ast},T_k),&\\
&(T_k^{\ast},T_k^{\ast})\xrightarrow{h_k}(U_{k1},T_k^{\ast})\xrightarrow{f_{k2}}(T_k,T_k^{\ast}).&
\end{eqnarray*}
\item[(IV)] For any $i\in Q_0\backslash F$ with $i\neq k$, we have that $T_k\notin(\add U_{i1})\cap(\add U_{i0})$ and the following sequences are exact.
\begin{eqnarray*}
&(T_k^{\ast},U_{i1})\xrightarrow{f_{i1}}(T_k^{\ast},U_{i0})\xrightarrow{f_{i0}}(T_k^{\ast},T_i),&\\
&(U_{i0},T_k^{\ast})\xrightarrow{f_{i1}}(U_{i1},T_k^{\ast})\xrightarrow{f_{i2}}(T_i,T_k^{\ast}).&
\end{eqnarray*}
\end{itemize}

\medskip
Now let $(Q',W') = \widetilde{\mu}_k(Q,W)$.
By Lemma \ref{ext_cont}, we can define a $\k$-linear functor $\Phi' \colon \widehat{\k Q'}\to\C$ in the following way.
\begin{itemize}
\item[(i)] $\Phi'$ coincides with $\Phi$ on $Q\cap Q'$.
\item[(ii)] Let $\Phi'[ab] =\Phi a\Phi b$ for each pair of arrows $a \colon i\to k$ and
$b \colon k\to j$ in $Q$.
\item[(iii)] Let $\Phi'k^{\ast}=T_k^{\ast}$, and $\Phi'a^{\ast}$ ($e(a)=k$) and $\Phi'b^{\ast}$ ($s(b)=k$) are defined by 
\begin{eqnarray*}
\Phi'((a^{\ast})_{a\in Q_1,\ e(a)=k})=g_k&\in&\C(T_k^{\ast},\bigoplus_{a\in Q_1,\ e(a)=k}T_{s(a)}),\\
\Phi'((b^{\ast})_{b\in Q_1,\ s(b)=k})=-h_k&\in&\C(\bigoplus_{b\in Q_1,\ s(b)=k}T_{e(b)},T_k^{\ast})
\end{eqnarray*}
which are given in (I).
\end{itemize}
As before we simply denote $\Phi'p$ by $p$ for any morphism $p$ in $\widehat{KQ'}$.

Our key result is the following:

\begin{theorem}\label{AR in T'}
For any vertex $i$ in $Q'_0\backslash F$, we have the following right 2-almost split sequence in $\add \mu_k(T)$.
\[T_i\xrightarrow{(d)_d}\bigoplus_{d\in Q'_1,\ s(d)=i}T_{e(d)}\xrightarrow{(\partial_{(c,d)}W')_{d,c}}
\bigoplus_{c\in Q'_1,\ e(c)=i}T_{s(c)}\xrightarrow{(c)_c}T_i.\]
\end{theorem}

Before proving Theorem \ref{AR in T'}, we notice that our main
Theorems \ref{main2}, \ref{main} and \ref{main3} follow immediately from Theorem \ref{AR in T'}.
We need the following observation.

\begin{lemma}
Let $\C$, $T$ and $(Q,W,F)$ be those in Theorems \ref{main2}, \ref{main} and \ref{main3},
where we change the notations as follows.
\[\begin{array}{|c||c|c|c|}\hline
     &\mbox{Theorem \ref{main2}}&\mbox{Theorem \ref{main}}&\mbox{Theorem \ref{main3}}\\ \hline \hline
\C   &\mod\Lambda               &\C                       &\E                         \\ \hline
(Q,W)&(Q,W)                     &(Q,W)                    &(\widetilde{Q},\widetilde{W})\\ \hline
F    &\emptyset                 &\emptyset                &F                           \\ \hline
\end{array}\]
If $(Q,W)$ is reduced, then the above conditions (O)--(IV) are satisfied.
\end{lemma}

\begin{proof}
(O) For the case Theorems \ref{main} and \ref{main2}, this follows from Theorem \ref{key criterion}.
For the case Theorem \ref{main3}, this is a part of the assumption.

(I) Since $(Q,W)$ is reduced, $f_{i1}\in J_{\C}$ holds.
Since $Q$ has no loops,
there exist exchange triangles/sequences
\[T_i\xrightarrow{f_{i2}}U_{i1}\xrightarrow{h_i}T_i^{\ast}\ \mbox{ and }\ T_i^{\ast}\xrightarrow{g_i}U_{i0}\xrightarrow{f_{i0}}T_i\]
such that $h_ig_i=f_{i1}$ by Theorem \ref{from 2-slmost split to exchange}.
In particular, we have (I).

(II) This follows from Proposition \ref{gluing of exchange}.

(III) This is clear from a property of triangles and exact sequences.

(IV) Since no 2-cycle starts at the vertex $k$, we have $T_k\notin(\add U_{i1})\cap(\add U_{i0})$.
Since $(T_k^{\ast},T_i^{\ast})\xrightarrow{g_i}(T_k^{\ast},U_{i0})\xrightarrow{f_{i0}}(T_k^{\ast},T_i)$ is exact, we only have to show that
$(T_k^{\ast},U_{i1})\xrightarrow{h_i}(T_k^{\ast},T_i^{\ast})$ is surjective.
This is clear for the case Theorems \ref{main} and \ref{main3} since we have $\Ext^1_{\C}(T_k^{\ast},T_i)=0$.

We consider the case Theorem \ref{main2}.
Fix any $p\in\Hom_\Lambda(T_k^{\ast},T_i^{\ast})$. Since $0\to(T,T_i)\xrightarrow{f_{i2}}(T,U_{i1})\xrightarrow{h_i}(T,T_i^{\ast})\to0$
is exact by Proposition \ref{3-CY mutation}(c),
there exist $q\in\Hom_\Lambda(U_{k1},U_{i1})$ and $r\in\Hom_\Lambda(T_k,T_i)$ which make the diagram
\[\xymatrix{
0\ar[r]&T_k\ar[r]^{f_{k2}}\ar[d]^r&U_{k1}\ar[r]^{h_k}\ar[d]^q&T_k^{\ast}\ar[d]^p\\
0\ar[r]&T_i\ar[r]^{f_{i2}}&U_{i1}\ar[r]^{h_i}&T_i^{\ast}
}\]
commutative. Since $0\to(T_k^{\ast},T)\xrightarrow{h_k}(U_{k1},T)\xrightarrow{f_{k2}}J_{\mod\Lambda}(T_k,T)\to0$ is exact by Proposition \ref{3-CY mutation}(c),
there exist $s\in\Hom_\Lambda(U_{k1},T_i)$ such that $r=f_{k2}s$
and $t\in\Hom_\Lambda(T_k^{\ast},U_{i1})$ such that $q=sf_{i2}+h_kt$.
Then we have $h_k(p-th_i)=0$.
By Proposition \ref{3-CY mutation}(b), we have $p=th_i$.
\end{proof}

Now we are ready to prove Theorems \ref{main2}, \ref{main} and \ref{main3}.
Without loss of generality, we can assume that $(Q,W)$ is reduced since $\mu_k(Q,W)$
is right equivalent to $\mu_k(Q_{\rm red},W_{\rm red})$.
Theorems \ref{main2} and \ref{main} follow from Theorems \ref{AR in T'} and \ref{key criterion}.
Theorem \ref{main3} follows immediately from Theorem \ref{AR in T'}.
\qed

\medskip
In the rest of this section, we prove Theorem \ref{AR in T'}.
We divide our proof of Theorem \ref{AR in T'} into Lemmas \ref{analogue of case 1}, \ref{analogue of case 2} and \ref{analogue of case 3}.
We need the following information about mutation of QP's.

\begin{lemma}\label{calculating partial}
Let $(Q,W)$ be a QP and $(Q',W') = \widetilde{\mu_k}(Q,W)$.
Let $a$ and $b$ be arrows in $Q$ with $e(a)=k=s(b)$, and let $c$ and $c'$ be arrows in $Q\cap Q'$.
Then we have the following equalities.
\begin{itemize}
\item[(a)] $\partial_{(c,c')}W'= \partial_{(c,c')}[W]$.
\item[(b)] $\partial_{(c, [ab])}W'=\partial_{(c, [ab])}[W]$ and $\partial_{([ab],c)}W'=\partial_{([ab],c)}[W]$.
\item[(c)] $\partial_{(a^{\ast}, [ab])}W'=b^{\ast}$.
\item[(d)] $\partial_{([ab],b^{\ast})}W'=a^{\ast}$.
\item[(e)] $\partial_{(b^{\ast}, a^{\ast})}W'=[ab]$.
\item[(f)] For other pairs $d,d' \in Q_1'$, we have $\partial_{(d,d')}W' = 0$.
\end{itemize}
\end{lemma}

\begin{proof}
Immediate from the definition $W'=[W]+\Delta$.
\end{proof}

\begin{lemma}\label{analogue of case 1}
We have the following weak 2-almost split sequence in $\add \mu_k(T)$.
\[T_k^{\ast}\xrightarrow{(a^{\ast})_a}{\bigoplus_{a\in Q_1,\ e(a)=k}T_{s(a)}\xrightarrow{([ab])_{a,b}}
\bigoplus_{b\in Q_1,\ s(b)=k}T_{e(b)}\xrightarrow{(b^{\ast})_b}T_k^{\ast}.}\]
\end{lemma}

\begin{proof}
By our definition of $\Phi'$, we can write the above sequence as
\[T_k^{\ast}\xrightarrow{g_k}U_{k0}\xrightarrow{f_{k0}f_{k2}}U_{k1}\xrightarrow{-h_k}T_k^{\ast}.\]
This is a weak 2-almost split sequence by (II).
\end{proof}

\begin{lemma}\label{analogue of case 2}
For a vertex $i$ in $Q'_0\backslash F$ with $i\neq k^{\ast}$, assume that there is no arrow $i\to k$ in $Q$.
Then we have a weak 2-almost split sequence
\begin{equation}\label{new 2-almost split}
T_i\xrightarrow{((b^{\ast})_b\ (d)_d)}{\displaystyle
\begin{array}{rr}
{\displaystyle(\bigoplus_{\begin{smallmatrix}b\in Q_1\\ b \colon k\to i\end{smallmatrix}}T_k^{\ast})}\\
{\displaystyle\oplus(\bigoplus_{\begin{smallmatrix}d\in Q_1\\ s(d)=i\end{smallmatrix}}T_{e(d)})}
\end{array}
\xrightarrow{e}
\begin{array}{rr}
{\displaystyle(\bigoplus_{\begin{smallmatrix}a,b\in Q_1\\ e(a)=k\\ b \colon k\to i\end{smallmatrix}}T_{s(a)})}\\
{\displaystyle\oplus(\bigoplus_{\begin{smallmatrix}c\in Q_1\\ s(c)\neq k\\ e(c)=i\end{smallmatrix}}T_{s(c)})}
\end{array}
\xrightarrow{{([ab])_{(a,b)}\choose (c)_c}}T_i}
\end{equation}
in $\add \mu_k(T)$, where
\begin{eqnarray*}
e=\left(\begin{array}{cc}(a^{\ast})_a&0\\ (\partial_{([ab],d)}[W])_{d,(a,b)}&(\partial_{(c,d)}[W])_{d,c}\end{array}\right).
\end{eqnarray*}
\end{lemma}

\begin{proof}
This is a complex since we have
\begin{eqnarray*}
&\sum_{\begin{smallmatrix}a\in Q_1\\ e(a)=k\end{smallmatrix}}a^*[ab]=\sum_{\begin{smallmatrix}a\in Q_1\\ e(a)=k\end{smallmatrix}}(a^*a)b=g_kf_{k0}b=0&\\
&\sum_{\begin{smallmatrix}a,b\in Q_1\\ e(a)=k\\ b \colon k\to i\end{smallmatrix}}(\partial_{([ab],d)}[W])[ab]+
\sum_{\begin{smallmatrix}c\in Q_1\\ s(c)\neq k\\ e(c)=i\end{smallmatrix}}(\partial_{(c,d)}[W])c
=\partial_dW=0,&\\
&(b^{\ast}a^{\ast}+\sum_{\begin{smallmatrix}d\in Q_1\\ s(d)=i\end{smallmatrix}}d(\partial_{([ab],d)}[W]))_{(a,b)}
=-h_kg_k+(\partial_{(a,b)}W)_{(a,b)}=-f_{k1}+f_{k1}=0,&\\
&\sum_{\begin{smallmatrix}d\in Q_1\\ e(d)=i\end{smallmatrix}}d(\partial_{(c,d)}[W])=\partial_cW=0.&
\end{eqnarray*}
To simplify notations, we decompose $U_{i0}=T_k^\ell\oplus U''_{i0}$ with $T_k\notin\add U''_{i0}$
and write
\begin{eqnarray*}
f_{i1}=(f'_{i1}\ f''_{i1})&:&U_{i1}\to U_{i0}=T_k^\ell\oplus U''_{i0},\\
f_{i0}={f'_{i0}\choose f''_{i0}}&:&U_{i0}=T_k^\ell\oplus U''_{i0}\to T_i.
\end{eqnarray*}
Then we have the following exact sequences by (III) and (IV).
\begin{equation}\label{zero}
(\mu_k(T),T_k^{\ast})\xrightarrow{g_k}(\mu_k(T),U_{k0})\xrightarrow{f_{k0}}(\mu_k(T),T_k),
\end{equation}
\begin{equation}\label{first}
(\mu_k(T),U_{i1})\xrightarrow{(f'_{i1}\ f''_{i1})}(\mu_k(T),T_k^\ell\oplus U''_{i0})
\xrightarrow{{f'_{i0}\choose f''_{i0}}}(\mu_k(T),T_i),
\end{equation}
We have the following exact sequence by (IV).
\begin{equation}\label{second}
(T_k^\ell\oplus U''_{i0},\mu_k(T))\xrightarrow{(f'_{i1}\ f''_{i1})}(U_{i1},\mu_k(T))\xrightarrow{f_{i2}}(T_i,\mu_k(T)).
\end{equation}
We can write the sequence \eqref{new 2-almost split} as
\[T_i\xrightarrow{(s\ f_{i2})}T_k^{{\ast}\ell}\oplus U_{i1}
\xrightarrow{{g_k^\ell\ 0\ \choose t\ f''_{i1}}}
U_{k0}^\ell\oplus U''_{i0}
\xrightarrow{{f_{k0}^\ell f'_{i0}\choose f''_{i0}}}T_i\]
where $tf_{k0}^\ell=f'_{i1}$ holds.

(i) We will show that ${f_{k0}^\ell f'_{i0}\choose f''_{i0}}$ is right almost split in $\add\mu_k(T)$.

First we will show that any morphism 
$p\in J_{\C}(T/T_k,T_i)$ factors through ${f_{k0}^\ell f'_{i0}\choose f''_{i0}}$. Since $f_{i0}={f'_{i0}\choose f''_{i0}}$ is right almost split in $\add T$,
there exists $(p_1\ p_2)\in\C(T/T_k,T_k^\ell\oplus U''_{i0})$ such that $p=p_1f'_{i0}+p_2f''_{i0}$.
Since $f_{k0}$ is right almost split in $\add T$, there exists $q\in\C(T/T_k,U_{k0}^\ell)$ such that $p_1=qf_{k0}^\ell$.
Then we have $p=(q\ p_2){f_{k0}^\ell f'_{i0}\choose f''_{i0}}$.
\[\xymatrix@C=2cm{
&&T/T_k\ar^p[dd]\ar_{p_2}[dl]\ar^{p_1}[ddl]\ar_q@/_3pc/[ddll]\\
&U''_{i0}\ar_{f''_{i0}}[dr]&\\
U_{k0}^\ell\ar_{f_{k0}^\ell}[r]&T_k^\ell\ar_{f'_{i0}}[r]&T_i
}\]

Next we take any $p\in J_{\C}(T_k^{\ast},T_i)$.
Since $g_k$ is left almost split in $\add\mu_k(T)$, there exists
$q\in\C(U_{k0},T_i)$ such that $p=g_kq$.
Since $T_k\notin\add U_{i1}$, we have $T_i\notin\add U_{k0}$.
Thus $q\in J_{\C}(U_{k0},T_i)$ holds, and the first case implies that
$q$ factors through ${f_{k0}^\ell f'_{i0}\choose f''_{i0}}$.
Hence $p$ factors through ${f_{k0}^\ell f'_{i0}\choose f''_{i0}}$.
\[\xymatrix@C=2cm{
T_k^{\ast}\ar_p[d]\ar^{g_k}[r]&U_{k0}\ar^q[dl]\\
T_i}\]

(ii) We will show that ${g_k^\ell\ 0\ \choose\ t\ f''_{i1}}$ is a pseudo-kernel of ${f_{k0}^\ell f'_{i0}\choose f''_{i0}}$ in $\add\mu_k(T)$.

Assume $(p_1\ p_2)\in\C(\mu_k(T),U_{k0}^\ell\oplus U''_{i0})$ satisfies $(p_1\ p_2){f_{k0}^\ell f'_{i0}\choose f''_{i0}}=0$.
Since $(p_1\ p_2){f_{k0}^\ell\ 0\choose\ 0\ \ 1}{f'_{i0}\choose f''_{i0}}=0$ holds and \eqref{first} is exact,
there exists $q\in\C(\mu_k(T),U_{i1})$ such that $q(f'_{i1}\
f''_{i1})=(p_1\ p_2){f_{k0}^\ell\ 0\choose\ 0\ \ 1}$.
Hence we have $qf'_{i1}=p_1f_{k0}^\ell$ and $qf''_{i1}=p_2$.
Since $(p_1-qt)f_{k0}^\ell=qf'_{i1}-qf'_{i1}=0$ and \eqref{zero} is exact,
there exists $r\in\C(\mu_k(T),T_k^{{\ast}\ell})$ such that $p_1-qt=rg_k^\ell$.
Then we have $(p_1\ p_2)=(r\ q){g_k^\ell\ 0\ \choose\ t\ f''_{i1}}$.
\[\xymatrix@C=2cm{
&\mu_k(T)\ar^0@/^5pc/[ddr]\ar^{p_1}[d]\ar^(0.3){p_2}@/^2pc/[dd]\ar_q@/_5pc/[ddl]\ar_r[dl]\\
T_k^{{\ast}\ell}\ar^{g_k^\ell}[r]&U_{k0}^\ell\ar^{f_{k0}^\ell}[r]&T_k^\ell\ar_{f'_{i0}}[d]\\
U_{i1}\ar_{f''_{i1}}[r]\ar^{f'_{i1}}[rru]\ar^t[ru]&U''_{i0}\ar_{f''_{i0}}[r]&T_i
}\]

(iii) We will show that the map
\[s=(b^{\ast})_b\colon (\C/J_{\C})(T_k^{{\ast}\ell},T_k^{\ast})\to (J_{\C}/J_{\add \mu_k(T)}^2)(T_i,T_k^{\ast})\]
is bijective.
By the assumption (C2), we have $\k=(\C/J_{\C})(T_k^{\ast},T_k^{\ast})$.
By (II), we have that $h_k=(b^{\ast})_b \colon \bigoplus_{b\in Q_1,\ s(b)=k}T_{e(b)} \to T_k^{\ast}$ is minimal right almost split in $\add \mu_k(T)$
since the middle morphism $f_{k0}f_{k2}$ in the sequence in (II) belongs to $J_{\C}$.
Thus we have that
$(J_{\C}/J_{\add \mu_k(T)}^2)(T_i,T_k^{\ast})$ is a $\k$-vector space with basis $\{b^{\ast} \mid b \colon k\to i \}$.
Thus the above map is bijective.

(iv) We will show that $(s\ f_{i2})$ is left almost split in $\add\mu_k(T)$.

Since $f_{i2}$ is left almost split in $\add T$, any morphism in
$J_{\C}(T_i,T/T_k)$ factors through $(s\ f_{i2})$. Then
take any $p\in J_{\C}(T_i,T_k^{\ast})$. By (iii), there
exists $p_1\in\C(T_k^{{\ast}\ell},T_k^{\ast})$ such that $p-sp_1\in J_{\add\mu_k(T)}^2(T_i,T_k^{\ast})$.
Since $h_k$ is right almost split in $\add\mu_k(T)$ by (II), there exists
$q\in J_{\C}(U_{k0},T_i)$ such that $p-sp_1=qh_k$.
Since $f_{i2}$ is left almost split in $\add T$, there exists $r\in\C(U_{i2},U_{k2})$ such that $q=f_{i2}r$.
Then we have $p=sp_1+f_{i2}rh_k=(s\ f_{i2}){p_1\choose rh_k}$.
\[\xymatrix@C=2cm{
U_{k1}\ar[r]^{h_k}&T_k^{\ast}&T_k^{{\ast}\ell}\ar_{p_1}[l]\\
&T_i\ar[ul]_q\ar[u]^p\ar[ur]^s\ar[r]^{f_{i2}}&U_{i2}\ar^r@/^3pc/[llu]
}\]

(v) We will show that ${g_k^\ell\ 0\ \choose t\ f''_{i1}}$ is a
pseudo-cokernel of $(s\ f_{i2})$ in $\add\mu_k(T)$.

Assume ${p_1\choose p_2}\in\C(T_k^{{\ast}\ell}\oplus U_{i1},T')$ with
$T'\in\add\mu_k(T)$ satisfies $(s\ f_{i2}){p_1\choose p_2}=0$.
We first show that there exists $q\in\C(U_{k0}^\ell,T')$ such that $p_1=g_k^\ell q$.
Since $g_k$ is left almost split in $\add\mu_k(T)$ by (II)$^{\op}$, we only have to show $p_1\in J_{\C}$.
We have to consider the case $T'=T_k^{\ast}$. Since $sp_1=-f_{i2}p_2\in J_{\add\mu_k(T)}^2$, we have $p_1\in J_{\C}$ by (i).

Since $f_{i2}(p_2-tq)=f_{i2}p_2+sg_k^\ell q=f_{i2}p_2+sp_1=0$ and \eqref{second} is exact,
there exists ${q_1\choose q_2}\in\C(T_k^\ell\oplus U''_{i0},T')$ such that $p_2-tq=(f'_{i1}\ f''_{i1}){q_1\choose q_2}$.
Then $p_2=tq+f'_{i1}q_1+f''_{i1}q_2=t(q+f_{k0}^\ell q_1)+f''_{i1}q_2$ holds.
Thus we have
${p_1\choose p_2}={g_k^\ell\ 0\ \choose t\ f''_{i1}}{q+f_{k0}^\ell q_1\choose q_2}$.
\[\xymatrix@C=2cm{
&T'&\\
T_k^{{\ast}\ell}\ar[r]^{g_k^\ell}\ar[ru]^{p_1}&U_{k0}^\ell\ar^{f_{k0}^\ell}[r]\ar_{q}[u]&T_k^\ell\ar_{q_1}[ul]\\
T_i\ar_{f_{i2}}[r]\ar_s[u]\ar^0@/^5pc/[uur]&U_{i1}\ar_{f''_{i1}}[r]\ar_t[u]\ar_{f'_{i1}}[ur]\ar^(0.3){p_2}@/^2pc/[uu]&U''_{i0}\ar_{q_2}@/_5pc/[uul]}
\]
\end{proof}

\begin{lemma}\label{analogue of case 3}
For a vertex $i$ in $Q'_0\backslash F$ with $i\neq k^{\ast}$, assume that there is no arrow $k\to i$ in $Q$.
Then we have a weak 2-almost split sequence
\[T_i\xrightarrow{(([ab])_{a,b}\ (d)_d)}{\displaystyle
\begin{array}{rr}{\displaystyle(\bigoplus_{\begin{smallmatrix}a,b\in Q_1\\ s(b)=k\\ 
a \colon i\to k\end{smallmatrix}}T_{e(b)})}\\
{\displaystyle\oplus(\bigoplus_{\begin{smallmatrix}d\in Q_1\\ s(d)=i\\ e(d)\neq k\end{smallmatrix}}T_{e(d)})}\end{array}
\xrightarrow{e}
\begin{array}{rr}{\displaystyle(\bigoplus_{\begin{smallmatrix}a\in Q_1\\ a \colon i\to k\end{smallmatrix}}T_k^{\ast})}\\
{\displaystyle\oplus(\bigoplus_{\begin{smallmatrix}c\in Q_1\\ e(c)=i\end{smallmatrix}}T_{s(c)})}
\end{array}
\xrightarrow{{(a^{\ast})_a\choose (c)_c}}T_i}\]
in $\add \mu_k(T)$, where
\begin{eqnarray*}
e=\left(\begin{array}{cc}(b^{\ast})_b&(\partial_{(c,[ab])}[W])_{(a,b),c}\\ 0&(\partial_{(c,d)}[W])_{d,c}\end{array}\right)
\end{eqnarray*}
\end{lemma}

\begin{proof}
This is a dual of Lemma \ref{analogue of case 2}.
\end{proof}

\section{2-CY-tilted algebras associated with elements in Coxeter groups}\label{sec5}

In this section we show that a large class of 2-CY-tilted algebras, including the cluster-tilted algebras
and a class of 2-CY-tilted algebras coming from stable categories of preprojective algebras
of Dynkin type, are given by QP's, by finding an explicit description of the potentials
for some of the 2-CY tilted algebras. We can then use Theorem \ref{main3} to get that all 2-CY-tilted algebras
in the same mutation
class are given by QP's. We prove that these QP's are rigid in the sense of \cite{dwz}.

These 2-CY-tilted algebras come from 2-CY triangulated categories constructed from elements
in the Coxeter groups associated with connected quivers with no loops, as investigated in
\cite{ir, birs}. (See \cite{gls2} for an alternative approach to the construction of a subclass 
of these 2-CY categories).
We start with recalling the relevant results from \cite{birs}, including a description of the 
quivers of some special 2-CY-tilted algebras.

Let $Q$ be a finite connected quiver with no loops, with
vertices $1, \dots ,n$ and set of arrows $Q_1$. The associated preprojective algebra
over the algebraically closed field $K$ is defined as follows.
For each arrow $a \in Q_1$ from $i$ to $j$, we add a corresponding arrow $a^{\ast}$ from
$j$ to $i$ to get a new quiver $\overline{Q}$. 
Then the \emph{preprojective algebra} is defined by
$$\Lambda=K\overline{Q}/\langle\sum_{a\in Q_1}(aa^{\ast}-a^{\ast}a)\rangle.$$
We shall write $(a^{\ast})^{\ast} = a$, and $\epsilon(a) = 1, \epsilon(a^{\ast}) = -1$.

The \emph{Coxeter group} $W_Q$ is presented by generators
$s_1,\dots,s_n$ with relations $s_is_j=s_js_i$ if there is no arrow
in $Q$ between $i$ and $j$ and $s_is_js_i=s_js_is_j$ if there is
precisely one arrow in $Q$ between $i$ and $j$.
Let $w$ be an element in the associated Coxeter group $W_Q$, and
$w=s_{u_1}\cdots s_{u_m}$ a reduced expression, where the $u_1, \dots, u_m$ are integers in  
$\{1, \dots, n\}$. For each integer $i$ in $\{1, \dots, n\}$, consider the ideal
$I_i = \Lambda (1 - e_i) \Lambda$ in $\Lambda$, where $e_i$ denotes the 
trivial path at the vertex $i$. Let $I_w = I_{u_1} \cdots I_{u_m}$, and let
$\Lambda_w = \Lambda/ I_w$. Then $I_w$ and $\Lambda_w$ are independent of the choice of reduced
expression for $w$, and $\Lambda_w$ is a finite dimensional $K$-algebra. Denote
by $\Sub \Lambda_w$ the full subcategory of $\mod \Lambda_w$ whose objects are
the submodules of finitely generated projective $\Lambda_w$-modules.
We have that the injective dimensions $\id_{\Lambda_w} \Lambda_w $ and $\id\Lambda_w{}_{\Lambda_w} $ are at most one, and hence
the stable category $\C = \underline{\Sub} \Lambda_w$ is a Hom-finite triangulated
2-CY category \cite{h}. Associated with each reduced expression 
$s_{u_1}\cdots s_{u_m}$ of $w$ is the object $T = \Lambda / I_{u_1} \oplus \Lambda / I_{u_1} I_{u_2}
\oplus \cdots \oplus \Lambda / I_{u_1} \cdots I_{u_m}$ in $\Sub \Lambda_w$.
Consider the algebra $\End_{\Lambda_w}(T)$ and the associated 2-CY-tilted algebra
$\End_{\C}(T)$. Alternatively we can write the basic part of $T$ as a sum of indecomposable
objects
$$\begin{array}{ccccccc}
T &=& T_1 &\oplus& T_2 &\oplus \cdots \oplus& T_m\\[2mm]
&=& \frac{P_{u_1}}{I_{u_1}P_{u_1}} 
&\oplus& \frac{P_{u_2}}{(I_{u_1}I_{u_2})P_{u_2}} &\oplus \dots \oplus&
\frac{P_{u_m}}{(I_{u_1} I_{u_2} \dots I_{u_m})P_{u_m}},
\end{array}$$
where $P_i$ is the indecomposable projective $\Lambda$-module associated with the vertex $i$.
The structure of $T_v$ is determined by the sequence of integers $u_1, \dots , u_v$.
Let $K_{v,1}$ be the smallest submodule of $P_{u_v}$
such that $P_{u_v}/K_{v,1}$ is a sum of copies of $S_{u_v}$ (in this case only $S_{u_v}$). Then let 
$K_{v,2}$ be the smallest submodule of $K_{v,1}$ such that $K_{v,1}/K_{v,2}$ is a sum of
copies of $S_{u_{v-1}}$ etc. Then $T_v=P_{u_v}/K_{v,v}$.

We also recall \cite{birs} the following description of the quiver $Q' = Q(u_1, \dots, u_m)$ of 
$\End_{\Lambda}(T)$. The vertices are $1,2,\cdots,m$ ordered from left to right.
We often denote the vertex $v$ by $i_r$ if it is the $r$-th vertex of type $i$. In this case we write $|i_r|=v$.
There are two kinds of arrows in $Q'$.
\begin{itemize}
\item[(i)] Arrows going to the left: For two consecutive vertices of type $i$, for $1 \leq i \leq n$,
draw an arrow from the right one to the left one.
\item[(ii)] Arrows going to the right: 
For each arrow $i \xrightarrow{a} j$ in $\overline{Q}$, draw an arrow $u \xrightarrow{a} v$ 
whenever the following are satisfied.
\begin{itemize}
\item[-] $u$ is of type $i$ and $v$ is of type $j$
\item[-] there is no vertex of type $i$ between $u$ and $v$
\item[-] if there is a vertex $v'$ of type $j$ after $v$, then there is a vertex of type $i$ between $v$ and $v'$
 \end{itemize}
\end{itemize}

The following picture gives an illustration.
$$\xymatrix@R0.5cm@C0.5cm{
&i_1&i_2\ar[l]&i_3\ar[l]\ar[rrrd]^a&&&&i_4\ar[llll]&i_5\ar[l]\ar[rd]^a&\\
j_1\ar[rrru]^{a^*}&&&&j_2\ar[llll]&j_3\ar[l]&j_4\ar[l]\ar[rru]^{a^*}&&&j_5\ar[lll]
}$$

\noindent{\bf Example }
We give a concrete example.
Let $Q$ be the quiver $\xymatrix@R0.1cm@C0.3cm{1\ar[r]^a&2\ar[r]^b&3\ar@/_1pc/[ll]_c}$,
and let $w=s_1s_2s_1s_3s_1s_2s_3s_1s_2s_3s_2$ be a reduced expression.
The associated object $T$ in $\Sub\Lambda_w$ is the following.
\begin{eqnarray*}
&\begin{smallmatrix}
1
\end{smallmatrix}
\oplus
\begin{smallmatrix}
2\\
&1
\end{smallmatrix}
\oplus
\begin{smallmatrix}
&1\\
2
\end{smallmatrix}
\oplus
\begin{smallmatrix}
&&3\\
&1&&2\\
2&&&&1
\end{smallmatrix}
\oplus
\begin{smallmatrix}
&&&1\\
&&2&&3\\
&&&1&&2\\
&&&&&&1
\end{smallmatrix}
\oplus
\begin{smallmatrix}
&&&2\\
&&3&&1\\
&1&&2&&3\\
2&&&&1&&2\\
&&&&&&&1
\end{smallmatrix}
\oplus
\begin{smallmatrix}
&&3\\
&1&&2\\
2&&3&&1\\
&1&&2&&3\\
&&&&1&&2\\
&&&&&&&1
\end{smallmatrix}
\oplus
\begin{smallmatrix}
&&&&1\\
&&&2&&3\\
&&3&&1&&2\\
&1&&2&&3&&1\\
2&&&&1&&2&&3\\
&&&&&&&1&&2\\
&&&&&&&&&&1
\end{smallmatrix}&\\
&\oplus
\begin{smallmatrix}
&&&2\\
&&3&&1\\
&1&&2&&3\\
2&&3&&1&&2\\
&1&&2&&3&&1\\
&&&&1&&2&&3\\
&&&&&&&1&&2\\
&&&&&&&&&&1
\end{smallmatrix}
\oplus
\begin{smallmatrix}
&&&&&3\\
&&&&1&&2\\
&&&2&&3&&1\\
&&3&&1&&2&&3\\
&1&&2&&3&&1&&2\\
2&&&&1&&2&&3&&1\\
&&&&&&&1&&2&&3\\
&&&&&&&&&&1&&2\\
&&&&&&&&&&&&&1
\end{smallmatrix}
\oplus
\begin{smallmatrix}
&&&&&&2\\
&&&&&3&&1\\
&&&&1&&2&&3\\
&&&2&&3&&1&&2\\
&&3&&1&&2&&3&&1\\
&1&&2&&&&1&&2&&3\\
2&&&&&&&&&&1&&2\\
&&&&&&&&&&&&&1
\end{smallmatrix}&
\end{eqnarray*}
The quiver $Q(1,2,1,3,1,2,3,1,2,3,2)$ is the following.
$$\xymatrix@R0.5cm@C0.5cm{
1_1\ar[rd]_{a_1}&&1_2\ar[ll]_{p_2}\ar[rdd]^{c^*_2}&&1_3\ar[ll]_{p_3}\ar[rd]^{a_3}\ar@/_1.5pc/[rrdd]^{c_3^*}
&&&1_4\ar[lll]_{p_4}\ar[rrrd]^{a_4}\ar@/_1.5pc/[rrdd]^{c^*_4}\\
&2_1\ar[rrru]^{a^*_1}\ar[rrd]_{b_1}&&&&2_2\ar[llll]^{q_2}\ar[rru]^{a^*_2}\ar[rd]^{b_2}&&&2_3\ar[lll]^(0.46){q_3}\ar[rd]^{b_3}&&2_4\ar[ll]^{q_4}\\
&&&3_1\ar[rru]^{b^*_1}\ar[ruu]^{c_1}&&&3_2\ar[lll]^{r_2}\ar[rru]^{b^*_2}\ar[ruu]^{c_2}&&&3_3\ar[lll]^{r_3}\ar[ru]_{b^*_3}
}$$

Recall that the corresponding maps are given as follows.

\begin{lemma}\label{arrow and map}
\begin{itemize}
\item[(a)] Let $i_r \overset{a}{\to} i_{r-1}$ be an arrow in $Q'$ going to the left.
Then the map $T_{i_r} \overset{a}{\to} T_{i_{r-1}}$ is given by the natural surjection.
\item[(b)] Let $i_r \overset{b}{\to} j_s$ be an arrow in $Q'$ corresponding to an arrow $i \overset{b}{\to} j$ in $\overline{Q}$.
Then the map $T_{i_r} \overset{b}{\to} T_{j_s}$ is given by multiplication with $b$.
\end{itemize}
\end{lemma}

For each $i$, the last (i.e., rightmost) vertex of type $i$ corresponds to the
indecomposable projective $\Lambda_w$-module 
$\frac{P_i}{I_wP_i}$. By dropping the last vertex of type
$i$ for each $i = 1, \dots, n$, where possible, we obtain the quiver $\underline{Q}' = \underline{Q}(u_1, \dots, u_m)$ 
of $\End_{\C}(T)$ from the quiver $Q'$ of $\End_{\Lambda}(T)$, and
$T$ is a cluster-tilting object in $\underline{\Sub} \Lambda_w$.





\medskip
Our aim is to show that the 2-CY-tilted algebra $\End_{\C}(T)$ is a Jacobian algebra $\P(\underline{Q}',W)$,
by giving an explicit description of the potential $W$. 
For each arrow $i_r \overset{b}{\to} j_s$ in $\underline{Q}'$,
we let $W_b= \epsilon(b) b b^{\ast} p$ if there is a (unique) arrow $j_s \overset{b^{\ast}}{\to} i_t$ in $\underline{Q}'$, where $p$ denotes the path $i_t \to i_{t-1} \to \cdots \to i_r$.
Otherwise we let $W_b=0$. Then let
$$W = \sum_{b\in \underline{Q}'_1} W_b.$$
Our strategy is to show that all the relations $\partial_a W$
are satisfied for $\End_{\C}(T)$, and then apply Proposition \ref{p:second} to show $\End_{\C}(T) \simeq \P(\underline{Q}',W)$.

\medskip\noindent{\bf Example }
Again we consider the above example. Deleting vertices $1_4$, $3_3$ and $2_4$, we have the quiver $\underline{Q}(1,2,1,3,1,2,3,1,2,3,2)$ as follows.
$$
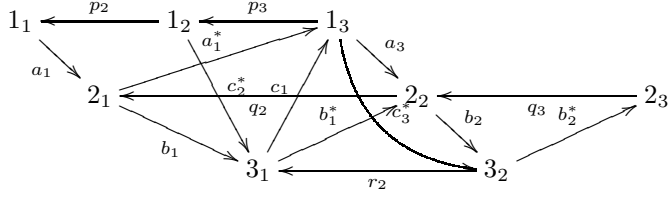
\xymatrix@R0.5cm@C0.5cm{
1_1\ar[rd]_{a_1}&&1_2\ar[ll]_{p_2}\ar[rdd]^{c^*_2}&&1_3\ar[ll]_{p_3}\ar[rd]^{a_3}\ar@/_1.5pc/[rrdd]^{c^*_3}&&&\ \ \ \\
&2_1\ar[rrru]^{a^*_1}\ar[rrd]_{b_1}&&&&2_2\ar[llll]^{q_2}\ar[rd]^{b_2}&&&2_3\ar[lll]^{q_3}&&\ \ \ \\
&&&3_1\ar[rru]^{b^*_1}\ar[ruu]^{c_1}&&&3_2\ar[lll]^{r_2}\ar[rru]^{b^*_2}&&&\ \ \ 
}$$
The associated potential is given by
\[W=a_1a_1^*p_3p_2-a_1^*a_3q_2+b_1b_1^*q_2-b_1^*b_2r_2+b_2b_2^*q_3-c_2^*c_1p_3+c_1c_3^*r_2.\]

\medskip
Let us start with the following information.

\begin{lemma}\label{arrow to the left}
Let $i_r\to i_{r-1}$ be an arrow in $Q'$ going to the left.
For any arrow $i \overset{b}{\to} j$ in $\overline{Q}$,
there exists a map $T_{i_{r-1}} \overset{bb^{\ast}}{\to} T_{i_r}$ of $\Lambda$-modules given by multiplication with $bb^{\ast}$.
Moreover precisely one of the following holds.
\begin{itemize}
\item[(a)] There do not exist vertices $i_u$ and $j_t$ in $Q'$ satisfying $|i_u|<|j_t|<|i_r|$.
In this case the map $T_{i_{r-1}} \overset{bb^{\ast}}{\to} T_{i_r}$ is zero.
\item[(b)] There exists a path $i_u \overset{b}{\to} j_t \overset{b^{\ast}}{\to} i_v$ in $Q'$ with $u\le r-1<r\le v$.
In this case such a path is unique, and the map $T_{i_{r-1}} \overset{bb^{\ast}}{\to} T_{i_r}$ is equal to the composition
$T_{i_{r-1}} \overset{p}{\to} T_{i_u} \overset{b}{\to} T_{j_t} \overset{b^{\ast}}{\to} T_{i_v} \overset{q}{\to} T_{i_{r}}$
where $p$ denotes a path $i_{r-1} \overset{}{\to} i_{r-2} \overset{}{\to} \cdots \overset{}{\to} i_u$
and $q$ denotes a path $i_v \overset{}{\to} i_{v-1} \overset{}{\to} \cdots \overset{}{\to} i_r$.
%
\end{itemize}
\end{lemma}

\begin{proof}
We can write $T_{i_{r-1}}=P_i/IP_i$ and $T_{i_r}=P_i/(II'I_i)P_i$, where $I'$ is a product of ideals $I_1,\ldots,I_n$ except $I_i$.
Then we have $bb^{\ast}\in I'I_i$. Thus we have $Ibb^{\ast}\subset II'I_i$, and the map $T_{i_{r-1}} \overset{bb^{\ast}}{\to} T_{i_r}$ is well-defined.

(a) Since $bb^{\ast}$ is zero in $T_{i_r}$, we have the assertion.

(b) The uniqueness of the path is clear from the definition of $Q'$.
The latter assertion is clear from Lemma \ref{arrow and map}.
\end{proof}

We now show that the 2-CY-tilted algebra $\End_{\C}(T)$ satisfies the relations for the Jacobian algebra 
$\P(\underline{Q}',W)$.

\begin{proposition}\label{p:third}
$\partial_a W$ belongs to the kernel of the surjection $\widehat{K \underline{Q}'}\to\End_{\C}(T)$ for any arrow $a$ in $\underline{Q}'$.
\end{proposition}

\begin{proof}
There are two cases to consider.

(1) Consider an arrow $i_r \overset{b}{\to} j_s$ going to the right in $\underline{Q}'$. If there is an arrow 
$j_s \overset{b^{\ast}}{\to} i_t$ in $\underline{Q}'$, we have the cycle 
$\epsilon(b)b b^{\ast} p$ as part of the potential $W$, where $p$ is the path $i_t \to i_{t-1} \to \cdots \to i_r$.
And if there is an arrow $j_u \overset{b^{\ast}}{\to} i_r$ in $\underline{Q}'$, we have the cycle
$\epsilon(b^{\ast})b^{\ast} b q$, where $q$ is the path $j_s \to j_{s-1} \to \cdots \to j_u$.
The relation $\partial_b W$ corresponding to the arrow $i_r \overset{b}{\to} j_s$ is one of the following four possibilities up to sign.
\begin{itemize}
\item[(i)] $b^{\ast} p - q b^{\ast}$, if both $j_s \overset{b^{\ast}}{\to} i_t$ and $j_u \overset{b^{\ast}}{\to} i_r$ exist,
\item[(ii)] $b^{\ast} p$, if $j_s \overset{b^{\ast}}{\to} i_t$ exists, but not $j_u \overset{b^{\ast}}{\to} i_r$,
\item[(iii)] $q b^{\ast}$, if $j_u \overset{b^{\ast}}{\to} i_r$ exists, but not $j_s \overset{b^{\ast}}{\to} i_t$,
\item[(iv)] $0$, if neither $j_s \overset{b^{\ast}}{\to} i_t$ nor $j_u \overset{b^{\ast}}{\to} i_r$ exist.
\end{itemize}
In each case, both $T_{j_s} \overset{b^{\ast} p}{\longrightarrow} T_{i_r}$ and $T_{j_s} \overset{q b^{\ast}}{\to} T_{i_r}$
are given by multiplication with $b^{\ast}$ by Lemma \ref{arrow and map}. Thus the cases (i) and (iv) are clear.
For the case (ii), there is no vertex $j_u$ of type $j$ satisfying $|j_u|<|i_r|$.
Since then $T_{i_r}$ has no composition factor $S_j$, we have that $T_{j_s} \overset{b^{\ast} p}{\to} T_{i_r}$ is zero.
For the case (iii), there exists an arrow $j_s \overset{b^{\ast}}{\to} i_t$ in $Q'$ such that $T_{i_t}$ is a projective $\Lambda_w$-module.
Thus $T_{j_s} \overset{q b^{\ast}}{\to} T_{i_r}$ is equal to $T_{j_s} \overset{b^{\ast} p}{\to} T_{i_r}$ which is zero in $\C$.

(2) Consider an arrow $i_{r} \overset{a}{\to} i_{r-1}$ going to the left in $\underline{Q}'$.

By definition of $W$ we have
$$\partial_aW=\sum_{b}\epsilon(b) pbb^{\ast}q,$$
where $b$ is an arrow in $\overline{Q}$ satisfying the condition in Lemma \ref{arrow to the left}(b) and the additional condition that
the path $pbb^{\ast}q$ is in $\underline{Q}'$.
In this case the map $T_{i_{r-1}} \overset{\epsilon(b) pbb^{\ast}q}{\longrightarrow} T_{i_r}$
is given by multiplication with $\epsilon(b)bb^{\ast}$ by Lemma \ref{arrow to the left}(b).

If $b$ is an arrow in $\overline{Q}$ satisfying the condition in Lemma \ref{arrow to the left}(b) and the path $pbb^{\ast}q$ is not in $\underline{Q}'$,
then the map $T_{i_{r-1}} \overset{\epsilon(b) pbb^{\ast}q}{\longrightarrow} T_{i_r}$ is zero in $\C$
and equal to the map which is multiplication with $\epsilon(b)bb^{\ast}$.

If $b$ is an arrow in $\overline{Q}$ which does not satisfy the condition in Lemma \ref{arrow to the left}(b),
then the multiplication map $T_{i_{r-1}} \overset{bb^{\ast}}{\to} T_{i_r}$ is zero by Lemma \ref{arrow to the left}(a).

Consequently the map $T_{i_{r-1}} \overset{\partial_aW}{\to} T_{i_r}$ is equal to multiplication with $\sum_{b\in\overline{Q}_1}\epsilon(b)bb^{\ast}$,
which is zero by the definition of the preprojective algebra.
\end{proof}


Now we are ready to prove one of the main results in this section.

\begin{theorem}\label{t:fourth}
Let $Q$ be a finite quiver without loops
and $W_Q$ the Coxeter group associated with $Q$. Let $\Lambda$ be the preprojective algebra
over the algebraically closed field $K$ and $w\in W_Q$. Let $T$ be a cluster-tilting
object associated with a reduced expression of $w$ in the 2-CY triangulated category
$\C = \underline{\Sub} \Lambda_w$.
Then the 2-CY-tilted algebra $\End_{\C}(T)$ is given by the QP $(\underline{Q}',W)$ above.
\end{theorem}

\begin{proof}
By Proposition \ref{p:third} we have that $\End_{\C}(T)$ is a factor algebra of the Jacobian algebra
$\P(\underline{Q}',W)$. We want to apply Proposition \ref{p:second}. 

We first claim that each cycle $C$ in $W$ is full. Recall 
that $C$ is of the form
$$
\xymatrix{
i_u \ar[drr] & \dots \ar[l] & \dots & i_{r-1} \ar[l] & i_r \ar[l] \\
& & j_s \ar[urr] & &
}
$$
There are no more arrows between the vertices of type $i$ in $C$.
Further $i_r$ is the unique vertex of type $i$ where there is an arrow from $j_s$, and $i_u$
is the unique vertex of type $i$ with an arrow to $j_s$. This shows that $C$ is full.

For any cycle $C$ of $W$ and any arrow $a$ in $C$, it is clear from the definition of $W$ that the path $\partial_aC$ determines $a$.
Hence for any distinct arrows $a$ and $b$ contained in cycles of $W$, we have that $\partial_aW$ and $\partial_bW$ do not have common paths.
Consequently the $\partial_aW$ for all arrows $a$ contained in cycles of $W$ are linearly independent over $K$.


Since $\Gamma = \End_{\C}(T)$ is 2-CY-tilted, it follows from \cite{kr1} that 
$\dim_K \Ext^2_{\Gamma}(S,S') \leq \dim_K \Ext^1_{\Gamma}(S',S)$ for simple $\End_{\C}(T)$-modules $S$ and $S'$.
It now follows from Proposition \ref{p:second} that $\End_{\C}(T)$ is isomorphic to the Jacobian algebra
$\P(\underline{Q}',W)$.
\end{proof}

We now show that the potential $W$ is rigid.

\begin{theorem}
With the previous notation, the potential $W$ on $\underline{Q}'$ is rigid.
\end{theorem}

\begin{proof}
Let $C$ be any cycle in the quiver $\underline{Q}'$. We want to show
that $C$ belongs to $\JJ(W)$ up to cyclic equivalences.

We say that a vertex $v$ in $C$ is a {\it right turning point}
if an arrow in $C$ going from left to right ends at $v$, and an arrow
in $C$ going from right to left starts at $v$. Then define
$r(C) = \sum_v 3^{|v|}$, where we sum over all right turning points 
$v$ of $C$.
Consider a vertex $i_u$ on $C$,
with $|i_u|$ minimal. Consider the last right turning point $i_v$ preceding $i_u$, and
let $j_s \to i_v$ be the preceding arrow in $C$. 
Since $|i_u| <|j_s| < |i_v|$, we can
choose $i_r$ with $|i_r| < |j_s| < |i_{r+1}|$.

Assume first that there is some $j_t$ with $|j_t| < |i_r|$. Then there is an 
arrow $j_t \overset{b}{\to} i_r$ if $j_t$ is chosen with $|j_t|$ largest possible.
We have the following subquiver of $\underline{Q}'$:
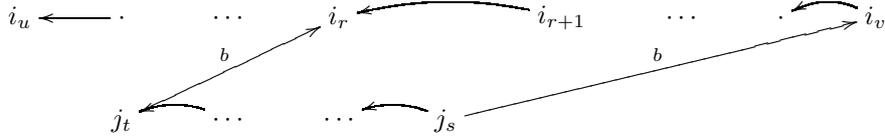
$$
\xymatrix{
i_u & \cdot \ar[l] & \cdots & i_r & & i_{r+1} \ar@/_/[ll] & \cdots & \cdot & i_v \ar@/_/[l] \\
& j_t \ar^{b}[urr] & \cdots \ar@/_/[l] & \cdots & j_s \ar^{b}[urrrr] \ar@/_/[l]  & & & &  
}
$$
The composition $j_s \overset{b}{\to} i_v \to i_{v-1} \to \cdots \to i_r$ coincides with 
the composition $j_s \to j_{s-1} \to \cdots \to j_t \overset{b}{\to} i_r$ by Lemma \ref{arrow and map}. 
So we replace the first path by the second one to get a new cycle $C'$ satisfying $C-C'\in\JJ(W)$.
The right turning point $i_v$ in $C$ is replaced by at most two 
right turning points $i_r$ and $j_s$ in $C'$, where $|i_r| < |i_v|$ and $|j_s| < |i_v|$,
and hence $3^{|i_r|} + 3^{|j_s|} < 3^{|i_v|}$. This shows that $r(C') < r(C)$,
and we are done by induction.

Assume now that there is no $j_t$ with $|j_t| < |i_r|$. Then $T_{i_r}$ has no
composition factor $S_j$, so the composition $T_{j_s} \to \cdots \to T_{i_r}$ must be 0.
Hence $C\in\JJ(W)$ in this case.
\end{proof}

Note that it follows from \cite{dwz} that all the potentials in the mutation class of 
$(\underline{Q}',W)$ are rigid, and hence we obtain a large class of
rigid QP's. Some examples of rigid QP's where the quiver $Q$ is not mutation
equivalent to an acyclic quiver were given in \cite{dwz}.

\medskip
There is a similar description of the (non-stable) endomorphism algebra $\End_\Lambda(T)$ for a cluster-tilting object $T$
in the 2-CY Frobenius category $\Sub\Lambda_w$ associated with a reduced expression of $w$,
using QP's with frozen vertices as discussed in section 1.

Let $w=s_{u_1}\cdots s_{u_m}$ be a reduced expression.
Let $T$ be the associated cluster-tilting object in $\Sub\Lambda_w$ and $Q'=Q(u_1,\ldots,u_m)$ the associated quiver.
We define a potential $W'$ of $Q'$ as follows.
For each arrow $i_r \overset{b}{\to} j_s$ in $Q'$,
we let $W'_b= \epsilon(b) b b^{\ast} p$ if there is a (unique) arrow $j_s \overset{b^{\ast}}{\to} i_t$ in $Q'$, where $p$ denotes the path $i_t \to i_{t-1} \to \cdots \to i_r$.
Otherwise we let $W'_b=0$. Then let
$$W' = \sum_{b\in Q'_1} W'_b.$$
Let $F$ be the set of vertices in $Q'$ which are not contained in $\underline{Q}'$. Then we have a QP $(Q',W',F)$ with frozen vertices.

We have the following analogue of Theorem \ref{t:fourth}.

\begin{theorem}\label{lifting of t:fourth}
$\End_{\Lambda}(T)$ is isomorphic to $\P(Q',W',F)$, and hence $\End_{\underline{\Sub}\Lambda_w}(T)$ is isomorphic to $\P(\underline{Q}',W)$.
\end{theorem}

We omit the proof since it is quite similar to that of Theorem \ref{t:fourth}.

Since the 2-CY categories $\underline{\Sub} \Lambda_w$ with the cluster-tilting objects
determine a cluster structure \cite{birs}, we have the following result.

\begin{cor}\label{mutation of reduced words are given by QP}
With the previous notation, all the 2-CY-tilted algebras
belonging to the mutation class of cluster-tilting objects associated with reduced expressions are given by rigid QP's. 
\end{cor}

\begin{proof}
The isomorphism $\P(Q,W)\to\End_\Lambda(T)$ is liftable by Theorems \ref{lifting of t:fourth}
and \ref{analogue of key criterion}.
Thus the assertion follows from Corollary \ref{c:mut-class}.
\end{proof}

Note that we also get a large class of liftable isomorphisms, which is useful for applying Corollary \ref{c:mut-class}.

Since the cluster categories of connected quivers are equivalent to $\Sub\Lambda_{c^2}$ 
for a Coxeter element $c$ of $Q$ (see \cite{birs,gls2}), and have connected cluster-tilting graph,
we get the following from Corollary \ref{mutation of reduced words are given by QP}.

\begin{cor}\label{jacob}
Any cluster-tilted algebra is isomorphic to a Jacobian algebra of a rigid QP.
 \end{cor}

This result can be used to give an alternative proof of the following.

\begin{cor}
Cluster-tilted algebras are determined by their quiver.
\end{cor}

\begin{proof}
If a cluster-tilted algebra $\Lambda$ has an acyclic 
quiver $Q$, then $\Lambda$ is isomorphic to the
path algebra $KQ$ by \cite{abs,kr1}, in particular it is determined by its
quiver.

Assume that $T_1$ and $T_2$ are cluster-tilting objects in
the cluster categories $\C_{Q_1}$ and $\C_{Q_2}$ respectively,
 such that the associated quivers are isomorphic. 
Let $\mu = \mu_{k_\ell}\circ \cdots \circ\mu_{k_1}$ be a sequence of
mutations such that for the projective $KQ_1$-module $KQ_1$ we have $\mu(KQ_1) \simeq T_1$,
and hence $\mu(Q_1) \simeq Q_{T_1}$. Then let $T$ be the cluster-tilting object 
in $\C_{Q_2}$ such that $\mu(T) \simeq T_2$ and hence $\mu(Q_T) \simeq Q_{T_2}$. Then
$Q_1 \simeq Q_T$, and hence $KQ_1 \simeq \End_{\C_{Q_2}}(T)$.
Since cluster-tilted algebras are QP's by Corollary \ref{jacob},
it follows from Theorem \ref{unique}(b) that $\End_{\C_1}(T_1) \simeq \End_{\C_2}(T_2)$.
\end{proof}

\noindent{\bf Remark }
It follows from the work in \cite{dwz} that if $\Lambda=\P(Q,W)$ is a Jacobian algebra 
given by a rigid potential $(Q,W)$ where $Q$ belongs to the mutation class of 
an acyclic quiver, then $\Lambda$ is determined by its quiver. Since by Corollary 
\ref{jacob} any cluster-tilted algebra has this property, this gives yet another 
way of seeing that cluster-tilted algebras are determined by their quiver.

\section{Nearly Morita equivalence for neighboring Jacobian algebras}\label{sec6}

Let $T$ be a cluster-tilting object in a triangulated 2-CY category $\C$, and $T^{\ast} = \mu_k(T)$
another cluster-tilting object obtained by mutation. Then we have a nearly Morita equivalence
between $\Lambda = \End_{\C}(T)$ and $\Lambda' = \End_{\C}(T^{\ast})$, that is an equivalence
$\frac{\mod \Lambda}{[\add S_k]} \to \frac{\mod \Lambda'}{[\add S_k']}$, 
where $S_k$ and $S_k'$ denote the simple modules at the vertex $k$ \cite{bmr1,kr1}.
This generalizes the equivalence of \cite{bgp} using reflection functors at sinks or sources.
Since there are Jacobian algebras which are not 2-CY-tilted, it is natural to ask if we have a nearly
Morita equivalence in general when performing mutations of quivers with potentials. It is the aim of this section 
to show that this is the case, when working with finite length modules.

Let $(Q,W)$ be a quiver with potential and $\Lambda=\P(Q,W)$ the associated Jacobian algebra.
Denote by $\widetilde{\mu}_k(Q,W)$ the quiver with potential and $\Lambda' = \P(\widetilde{\mu}_k(Q,W))$ the 
Jacobian algebra obtained by mutation at 
the vertex $k$. We show that there is an equivalence of categories
$\frac{\fd \Lambda}{[\add S_k]} \to \frac{\fd \Lambda'}{[\add S_k']}$, where $S_k$ and $S_k'$ denote the simple modules at the vertex $k$.

Our starting point is the map $G$ from objects in $\fd \Lambda$ to objects in $\fd \Lambda'$ used in \cite{dwz}, which we now recall.

Given a $\Lambda$-module $M$ and a vertex $k$, we let $a_1, ..., a_s$ be all arrows in $Q$ with $e(a_p)=k$ and $b_1, ..., b_t$ 
be all arrows with $s(b_q)=k$. We write
$$\Min=\bigoplus^s_{p=1}M_{e(a_p)}, \qquad \Mout=\bigoplus^t_{q=1}M_{s(b_q)}.$$

Let $\alpha \colon \Min \rightarrow M_k$ and $\beta \colon M_k\rightarrow \Mout$ be the maps given in matrix form by 
$$\alpha=(a_1, ..., a_s), \qquad \beta=\left(\begin{array}{c} b_1 \\ \vdots \\ b_t \end{array}\right).$$

Also, using the potential $W$, we define a map $\gamma \colon \Mout \rightarrow \Min$
$$\gamma=\left(\begin{array}{ccc} \gamma_{1,1} & \cdots & \gamma_{1,t}\\
\vdots & \ddots & \vdots\\
\gamma_{s,1} & \cdots & \gamma_{s,t}\end{array}\right)$$
where $\gamma_{p,q} = \partial_{(a_p,b_q)}W$ in the notation of section \ref{sec3}. 
Thus, locally, we get a triangle
$$\xymatrix{& M_k \ar[dr]^{\beta}\\
\Min \ar[ur]^{\alpha} && \Mout . \ar[ll]^{\gamma}}$$

Note that $\gamma_{p,q}$ is a linear combination of paths in $Q$ from $e(b_q)$ to $s(a_p)$, and that
there may be other paths from $\Mout$ to $\Min$.

Now, starting with $M$, we define a representation  $\widetilde{M}$ of $\Lambda' =\P(\widetilde{\mu}_k(Q,W))$ as follows: 
first, we let $i \colon \Ker \gamma \to \Mout$ be the inclusion and
$\varphi \colon \Ker \alpha \to \frac{\Ker \alpha}{\Im \gamma}$ the natural surjection, and
choose two $K$-linear splittings, 

$\rho \colon \Mout \rightarrow \Ker \gamma$ such that $i\rho = \id_{\Ker \gamma}$, and

$\sigma \colon \frac{\Ker \alpha}{\Im \gamma} \rightarrow \Ker \alpha$ such that $\sigma\varphi=\id_{\frac{\Ker \alpha}{\Im \gamma}}$.

Then, locally, the representation $\widetilde{M}$ is given by:
$$\xymatrix{& \widetilde{M_k} \ar[dl]_{\widetilde{\beta}}\\
\Min \ar@<.5ex>[rr]^{\alpha\beta} && \Mout \ar@<.5ex>[ll]^{\gamma} \ar[ul]_{\widetilde{\alpha}}}$$ 
where
$$\widetilde{M_k}=\frac{\Ker \gamma}{\Im \beta} \oplus \Im \gamma \oplus \frac{\Ker \alpha}{\Im \gamma} 
 \text{ , } \quad
 \widetilde{\alpha}=(-\rho\pi, -\gamma, 0) \quad \text{and} \quad \widetilde{\beta}=\left(\begin{array}{c}0\\ \iota\\ \sigma \end{array}\right),$$
where $\iota \colon \Im \gamma \rightarrow \Ker \alpha$ is the natural injection and
$\pi \colon \Ker\gamma \rightarrow \frac{\Ker\gamma}{\Im \beta}$ is the natural projection.

By \cite[Proposition 10.7]{dwz}, this makes $\widetilde{M}$ a representation of 
$\Lambda'$. Moreover, by \cite[Proposition 10.9]{dwz}, the isomorphism class of $\widetilde{M}$ 
does not depend on the choice of the splittings $\rho$ and $\sigma$. 

We next want to define an associated map on morphisms. 
Assume that we have two representations $M$ and $M'$ of $\Lambda$, and a map $f \colon M\rightarrow M'$, that is locally 
a commutative diagram 
\begin{equation}\label{commdia}
\xymatrix{&&& M_k \ar[dr]^{\beta} \ar[dd]_(.3){f_k}\\
M \ar[dd]_f && \Min \ar[ur]^{\alpha} \ar[dd]_{\fin} && \Mout \ar[ll]^(.3){\gamma} \ar[dd]^{\fout}
\\ &&& M'_k \ar[dr]^{\beta'} \\
M' && \Min' \ar[ur]^{\alpha'} && \Mout' . \ar[ll]^{\gamma'}
}
\end{equation}
We want to define a map $\widetilde{M}\rightarrow \widetilde{M'}$. 
To do so, we need to define a map $\widetilde{f_k} \colon \widetilde{M_k}\rightarrow \widetilde{M'_k}$ such 
that the squares of the following diagram commute: 
$$\xymatrix{&&& \widetilde{M_k} \ar[dl]_{\widetilde{\beta}} \ar[dd]_(.3){\widetilde{f_k}}\\
\widetilde{M} \ar[dd]_{\widetilde{f}} && \Min \ar@<.5ex>[rr]^(.7){\alpha\beta} \ar[dd]_{\fin} && \Mout \ar@<.5ex>[ll]^(.3){\gamma} 
\ar[dd]^{\fout} \ar[ul]_{\widetilde{\alpha}}
\\ &&& \widetilde{M'_k} \ar[dl]_{\widetilde{\beta'}} \\
\widetilde{M'} && \Min' \ar@<.5ex>[rr]^{\alpha'\beta'} && \Mout' \ar@<.5ex>[ll]^{\gamma'} \ar[ul]_{\widetilde{\alpha'}}
}$$
\noindent {\bf Remark }
>From the first diagram we have induced maps $\Ker \gamma \to \Ker \gamma'$, $\Ker \alpha \to \Ker \alpha'$,
$\Im \gamma \to \Im \gamma'$,
$\Im \beta \to \Im \beta'$, $\frac{\Ker \gamma}{\Im \beta} \to \frac{\Ker \gamma'}{\Im \beta'}$ and
$\frac{\Ker \alpha}{\Im \gamma} \to \frac{\Ker \alpha'}{\Im \gamma'}$.
This again induces componentwise a natural map $\widetilde{M_k} \to \widetilde{M_k'}$. However, this map
turns out not to work for our purposes, as the following example shows.\\
\ \\
\noindent {\bf Example }
Consider the Jacobian algebra $\P(Q,W)$ given by the quiver
$$
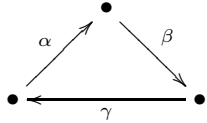
\xymatrix{& \bullet \ar[dr]^{\beta} & \\ 
\bullet \ar[ur]^{\alpha} & & \bullet \ar[ll]^{\gamma}
}$$
and the potential $W = \alpha \beta \gamma$.
Then $\widetilde{\mu}_k(\P(Q,W))= (\widetilde{Q}, \widetilde{W})$, where 
$\widetilde{Q}$ is the quiver 
$$
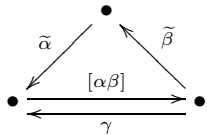
\xymatrix{& \bullet \ar[dl]_{\widetilde{\alpha}} & \\ 
\bullet \ar@<0.3ex>[rr]^{[\alpha \beta]} & & \bullet \ar@<1ex>[ll]^{\gamma} \ar[ul]_{\widetilde{\beta}}
}$$
and $\widetilde{W}$ is the the potential $[\alpha \beta] \gamma + [\alpha \beta] \widetilde{\alpha} \widetilde{\beta}$.
Consider a nonzero map between $\P(Q,W)$-modules
$$\xymatrix{& 0 \ar[dr] &  & \ar@<-4ex>[r] & & & 0 \ar[dr] & \\ 
K \ar[ur] & & K \ar[ll]^{\id} & & & 0 \ar[ur] & & K \ar[ll]  
}$$
The image of this map 
in $\fd \P(\widetilde{Q},\widetilde{W})$ would be 
$$\xymatrix{& (0,K,0) \ar[dl] &  & \ar@<-4ex>^{h}[r] & & & (0,0,K) \ar[dl] & \\ 
K \ar@<0.3ex>[rr]^{0} & & K \ar@<1ex>[ll]^{\id} \ar[ul] & & & 0 \ar@<0.3ex>[rr] & & K \ar@<1ex>[ll] \ar[ul]  
}$$
where $h$ is the zero-map, which indicates that this is not the map we are looking for.\\
\\ \
In view of the above, we now proceed to define a new map.
Considering the expression $$\frac{\Ker \gamma}{\Im \beta} \oplus \Im \gamma \oplus \frac{\Ker \alpha}{\Im \gamma},$$
we observe that the two first terms extend to $\Cok \beta$, while the two last terms 
extend to $\Ker \alpha$.  
We have a natural induced map $\Cok \beta \to \Ker \alpha$.
Taking these features into account we consider the  
%
%
knitting of short exact sequences
$$\xymatrix@R=15pt@C=15pt{0 \ar[rr] && \frac{\Ker \gamma}{\Im \beta} \ar[rr]^{\widetilde{i}} 
&& \Cok\beta \ar[rr]^{\widetilde{\gamma}\iota} \ar[dr]^{\widetilde{\gamma}} && \Ker \alpha \ar[rr]^{\varphi} && 
\frac{\Ker \alpha}{\Im \gamma} \ar[rr] && 0 \\
&&&&& \Im \gamma \ar[ur]^{\iota} \ar[dr]\\ 
&&&& 0\ar[ur] && 0}$$
Moreover, $\widetilde{i} \colon \frac{\Ker \gamma}{\Im \beta} \to \Cok \beta$ and
$\widetilde{\gamma} \colon \Cok \beta \to \Im \gamma$
are the induced maps.

Considering the induced map $\widetilde{\rho} \colon \Cok \beta \rightarrow \frac{\Ker\gamma}{\Im \beta}$, 
we get a knitting of ($K$-split) short exact sequences
$$\xymatrix@R=15pt@C=15pt{0 \ar[rr] && \frac{\Ker \gamma}{\Im \beta} \ar@<.5ex>[rr]^{\widetilde{i}} && \Cok\beta 
\ar@<.5ex>[ll]^{\widetilde{\rho}} \ar[rr]^{\widetilde{\gamma}\iota} \ar@<.5ex>[dr]^{\widetilde{\gamma}} && \Ker \alpha \ar@<.5ex>[rr]^{\varphi} 
\ar@<.5ex>[dl]^{\varepsilon} && \frac{\Ker \alpha}{\Im \gamma} \ar@<.5ex>[ll]^{\sigma} \ar[rr] && 0 \\
&&&&& \Im \gamma \ar@<.5ex>[ur]^{\iota} \ar@<.5ex>[ul]^{j} \ar[dr]\\ 
&&&& 0\ar[ur] && 0}$$
where the maps $j$ and $\varepsilon$ are induced from the splittings $\widetilde{\rho}$ and $\sigma$ respectively,
so that $\id_{\Cok \beta} = \widetilde{\rho}  \widetilde{i} +  \widetilde{\gamma} j$ and
$\id_{\Ker \alpha} = \varepsilon \iota + \varphi \sigma$.

The idea is then to use the diagram
$$\xymatrix@R=7pt{\frac{\Ker\alpha}{\Im \gamma} \ar[dr]^{\sigma} &&& \frac{\Ker\alpha'}{\Im \gamma'}\\
\oplus & \Ker\alpha \ar[r]^{\fin} & \Ker\alpha' \ar[ur]^{\varphi'}\ar[dr]^{\varepsilon'} & \oplus\\
\Im\gamma \ar[ur]^{\iota} \ar[dr]_j &&&\Im\gamma'\\
\oplus & \Cok\beta \ar[r]_{\widetilde{\fout}} \ar@{-->}[uu]_{\widetilde{\gamma}}& 
\Cok \beta' \ar@{-->}[uu]^{\widetilde{\gamma'}} \ar[ur]_{\widetilde{\gamma'}} \ar[dr]_{\widetilde{\rho'}}& \oplus\\
\frac{\Ker \gamma}{\Im \beta} \ar[ur]_{\widetilde{i}} &&& \frac{\Ker \gamma'}{\Im \beta'}}$$
where $\widetilde{\fout}$ is the induced quotient map and where the middle part commutes, that is 
$$\iota\circ\fin\circ\varepsilon'=j\circ\widetilde{\fout}\circ\widetilde{\gamma'}.$$ 

We let
\begin{multline*} 
\widetilde{f_k} = \left(%
\begin{array}{ccc} 
\widetilde{i}\circ\widetilde{\fout}\circ \widetilde{\rho'} & \widetilde{i}\circ \widetilde{\fout}\circ \widetilde{\gamma'} & 0 \\
j\circ\widetilde{\fout}\circ \widetilde{\rho'} & \iota\circ \fin \circ \varepsilon' & \iota\circ \fin \circ \varphi'\\
0 & \sigma\circ \fin \circ \varepsilon' & \sigma\circ \fin \circ \varphi'
\end{array}
\right)=  \\ \left(%
\begin{array}{ccc} 
\widetilde{i}\circ\widetilde{\fout}\circ \widetilde{\rho'} & 0 & 0 \\
j\circ\widetilde{\fout}\circ \widetilde{\rho'} & \iota\circ \fin \circ \varepsilon' & 0\\
0 & \sigma\circ \fin \circ \varepsilon' & \sigma\circ \fin \circ \varphi'
\end{array}
\right) \end{multline*}
since $\iota\circ \fin \circ \varphi'=0$ and $\widetilde{i}\circ \widetilde{\fout}\circ \widetilde{\gamma'}=0$. 
It is however often convenient to use the first expression of $\widetilde{f_k}$.

We need to verify that $\widetilde{\beta}\circ\fin  = \widetilde{f_k}\circ\widetilde{\beta'}$ and 
$\widetilde{\alpha}\circ\widetilde{f_k} = \fout\circ \widetilde{\alpha'}$.

First, 
$$ 
\begin{array}{rcl}
\widetilde{f_k}\circ\widetilde{\beta'} & = & \left(%
\begin{array}{ccc} 
\widetilde{i}\circ\widetilde{\fout}\circ \widetilde{\rho'} & \widetilde{i}\circ \widetilde{\fout}\circ \widetilde{\gamma'} & 0 \\
j\circ\widetilde{\fout}\circ \widetilde{\rho'} & \iota\circ \fin \circ \varepsilon' & \iota\circ \fin \circ \varphi'\\
0 & \sigma\circ \fin \circ \varepsilon' & \sigma\circ \fin \circ \varphi'
\end{array}
\right)\left(\begin{array}{c}0\\ \iota'\\ \sigma' \end{array}\right)\\ [1mm]
& = & \left(\begin{array}{c} 0 \\ \iota\circ\fin\circ(\varepsilon'\circ\iota' + \varphi'\circ\sigma') \\ [1mm]
\sigma\circ\fin\circ(\varepsilon'\circ\iota' + \varphi'\circ\sigma')
\end{array}\right) \\ [1mm]
& = & \left(\begin{array}{c}0\\ \iota\circ\fin\\ \sigma\circ\fin \end{array}\right) \\ [1mm]
& = & \widetilde{\beta}\circ\fin \end{array}
$$
%
On the one hand, we have
$$\fout\circ\widetilde{\alpha'}=(-\fout\circ\rho'\circ\pi', -\fout\circ\gamma', 0).$$
On the other hand, 
$$ 
\begin{array}{rcl}
\widetilde{\alpha}\circ\widetilde{f_k} & = & (- \rho\pi, -\gamma, 0)\left(%
\begin{array}{ccc} 
\widetilde{i}\circ\widetilde{\fout}\circ \widetilde{\rho'} & 0 & 0 \\
j\circ\widetilde{\fout}\circ \widetilde{\rho'} & \iota\circ \fin \circ \varepsilon' & 0\\
0 & \sigma\circ \fin \circ \varepsilon' & \sigma\circ \fin \circ \varphi'
\end{array}
\right)\\
& = & (-\rho\circ\pi\circ\widetilde{i}\circ\widetilde{\fout}\circ\widetilde{\rho'}-\gamma\circ j \circ\widetilde{\fout}\circ 
\widetilde{\rho'}, -\gamma\circ\iota\circ\fin\circ\varepsilon', 0)
\end{array}
$$
Since $-\gamma\circ\iota\circ\fin\circ\varepsilon'= -\gamma\circ\fin\circ\iota'\circ\varepsilon'=
-\gamma\circ\fin=-\fout\circ\gamma'$, it suffices (to get $\widetilde{\alpha}\circ\widetilde{f_k} = \fout\circ \widetilde{\alpha'}$) 
to show that 
$$-\rho\circ\pi\circ\widetilde{i}\circ\widetilde{\fout}\circ\widetilde{\rho'}-\gamma\circ j \circ\widetilde{\fout}\circ 
\widetilde{\rho'}=-\fout\circ\rho'\circ\pi'.$$
But, by the definition of $\widetilde{\rho}$ and $\widetilde{\gamma}$, the diagrams
$$\begin{array}{ccc} \xymatrix{\Mout\ar[r]^\rho \ar[d]_\pi & \Ker\gamma\ar[d]^\pi\\
\Cok\beta\ar[r]^{\widetilde{\rho}} & \frac{\Ker \gamma}{\Im\beta}} 
& \text{ and }& %
\xymatrix{\Mout \ar[r]^\gamma \ar[d]_{\pi} & \Im\gamma \\ \Cok\beta\ar[ur]_{\widetilde{\gamma}}} 
\end{array}$$
are commutative (where, by abuse of notation, we also denote by $\pi$ the natural projection $\Mout \rightarrow \Cok \beta$).
Therefore, 
$$\begin{array}{rcl}
-\rho\circ\pi\circ\widetilde{i}\circ\widetilde{\fout}\circ\widetilde{\rho'}-\gamma\circ j \circ\widetilde{\fout}\circ \widetilde{\rho'} & 
= & -(\rho\circ\pi\circ\widetilde{i} + \gamma\circ j)\circ \widetilde{\fout}\circ \widetilde{\rho'}\\
&=& -(\pi\circ\widetilde{\rho}\circ\widetilde{i} + \pi\circ\widetilde{\gamma}\circ j)\circ \widetilde{\fout}\circ \widetilde{\rho'}\\
&=& -\pi\circ(\widetilde{\rho}\circ\widetilde{i} + \widetilde{\gamma}\circ j)\circ \widetilde{\fout}\circ \widetilde{\rho'}\\
&=& -\pi\circ\id_{\Cok\beta}\circ \widetilde{\fout}\circ \widetilde{\rho'}\\
&=& -\pi\circ \widetilde{\fout}\circ \widetilde{\rho'}\\
&=& -\fout\circ\rho'\circ\pi'
\end{array}$$
where the last equality follows from the commutativity of 
$$\xymatrix{\Mout\ar[r]^{\fout}\ar[d]_\pi & \Mout'\ar[r]^{\rho'}\ar[d]^{\pi'} & \Ker\gamma' \ar[d]^{\pi'}\\
\Cok\beta\ar[r]^{\widetilde{\fout}} & \Cok\beta'\ar[r]^{\widetilde{\rho'}}&\frac{\Ker\gamma'}{\Im\beta'}}$$


Note that the image $\widetilde{f}$ apparently depends on the choice of splittings $\rho, \sigma, \rho'$
and $\sigma'$. However, we only aim to show that we obtain a functor
$\fd \Lambda \to \frac{\fd \Lambda'}{[\add S_k']}$. 

In \cite{dwz} it is proved that different choices of splittings will give isomorphic objects, and for each 
vertex $i \neq k$, the isomorphism is given by $\id \colon M_i \to \widetilde{M}_i$.
Therefore, different choices of splittings will give maps which are equal in $\frac{\fd \Lambda'}{[\add S_k']}$.
So, we get a well defined map of morphisms from $\fd \Lambda$ to 
$\frac{\fd \Lambda'}{[\add S_k']}$.

We next want to show that the definition of $G$ on objects and morphisms from $\fd \Lambda$ to $\frac{\fd \Lambda'}{[\add S_k']}$
actually gives rise to a functor, and that we get an induced equivalence
$$\frac{\fd \Lambda}{[\add S_k]} \to \frac{\fd \Lambda'}{[\add S_k']}.$$

We first observe that for $M$ in $\fd \Lambda$ we have $F(\id_M) = \id_{FM}$. 
This follows from $\sigma \varphi = \id_{\frac{\Ker \alpha}{\Im \gamma}}$,
$\widetilde{i} \widetilde{\rho} = \id_{\frac{\Ker \gamma}{\Im \beta}}$ and $\iota \varepsilon = \id_{\Im \gamma}$. Then, 
let $f \colon M\rightarrow M'$ and $f' \colon M'\rightarrow M''$ be two maps in $\fd \Lambda$.  
We will compare $\widetilde{f\circ f'}$ and $\widetilde{f}\circ\widetilde{f'}$. 
By construction, these maps are obtained from the following diagrams:
$$
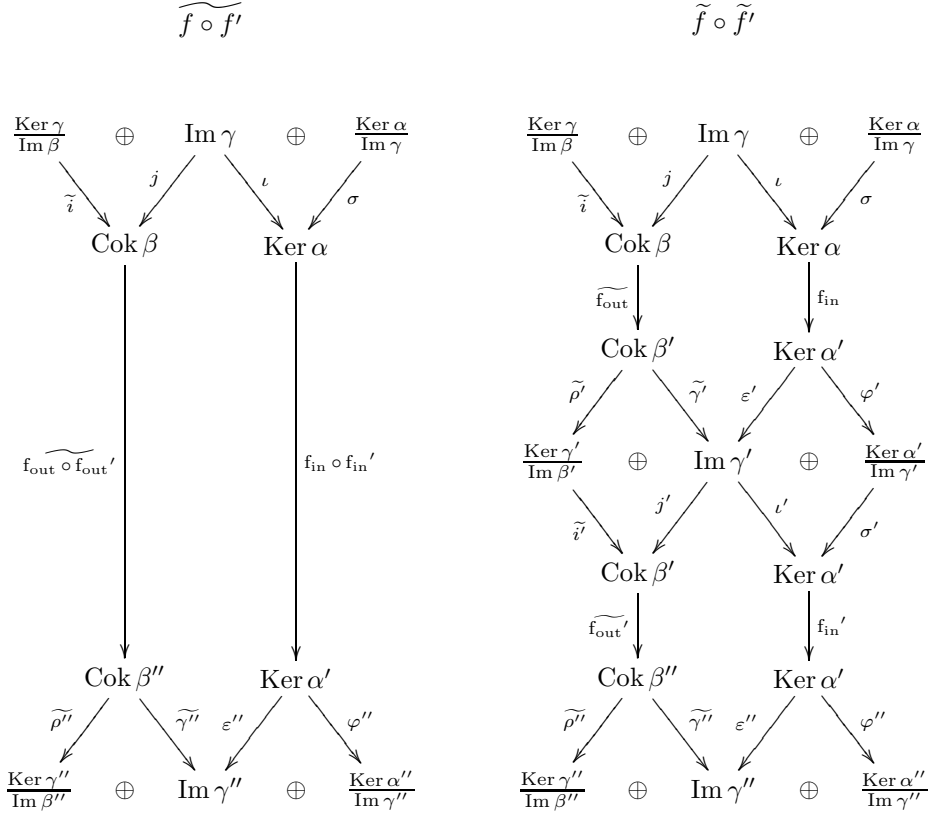
\xymatrix@!C=8pt{&&\widetilde{f\circ f'} &&&&&& \widetilde{f}\circ\widetilde{f'}\\
\frac{\Ker\gamma}{\Im\beta} \ar[dr]_{\widetilde{i}}  & \oplus & \Im\gamma\ar[dl]_j \ar[dr]^{\iota} & 
\oplus & \frac{\Ker\alpha}{\Im\gamma}\ar[dl]^\sigma && \frac{\Ker\gamma}{\Im\beta} \ar[dr]_{\widetilde{i}} & 
\oplus & \Im\gamma\ar[dl]_j \ar[dr]^{\iota} & \oplus & \frac{\Ker\alpha}{\Im\gamma}\ar[dl]^\sigma\\
&\Cok\beta\ar[dddd]_{\widetilde{\fout\circ\fout'}} &&\Ker\alpha\ar[dddd]^{\fin\circ\fin'}&&&& 
\Cok\beta\ar[d]_{\widetilde{\fout}} &&\Ker\alpha\ar[d]^{\fin}\\
&&&&&&&\Cok\beta'\ar[dl]_{\widetilde{\rho'}}\ar[dr]^{\widetilde{\gamma'}}&&\Ker\alpha'\ar[dl]_{\varepsilon'}\ar[dr]^{\varphi'} \\
&&&&&& \frac{\Ker\gamma'}{\Im\beta'} \ar[dr]_{\widetilde{i'}}  & \oplus & \Im\gamma'\ar[dl]_{j'} \ar[dr]^{\iota'} & 
\oplus & \frac{\Ker\alpha'}{\Im\gamma'}\ar[dl]^{\sigma'}\\
&&&&&&& \Cok\beta'\ar[d]_{\widetilde{\fout'}}&&\Ker\alpha'\ar[d]^{\fin'}\\
& \Cok\beta''\ar[dl]_{\widetilde{\rho''}}\ar[dr]^{\widetilde{\gamma''}}&&\Ker\alpha'\ar[dl]_{\varepsilon''}\ar[dr]^{\varphi''} &&&& 
\Cok\beta''\ar[dl]_{\widetilde{\rho''}}\ar[dr]^{\widetilde{\gamma''}}&&\Ker\alpha'\ar[dl]_{\varepsilon''}\ar[dr]^{\varphi''}\\
\frac{\Ker\gamma''}{\Im\beta''}  & \oplus & \Im\gamma''& \oplus & \frac{\Ker\alpha''}{\Im\gamma''} && \frac{\Ker\gamma''}{\Im\beta''} & 
\oplus & \Im\gamma''& \oplus & \frac{\Ker\alpha''}{\Im\gamma''}}
$$
and easy computations show that, locally, the morphism $\widetilde{f}\circ\widetilde{f'}-\widetilde{f\circ f'}$ is given by
$$\xymatrix@C=20pt@R=16pt{&\frac{\Ker\gamma}{\Im\beta}\oplus\Im\gamma\oplus\frac{\Ker\alpha}{\Im\gamma}\ar[dl]_{\widetilde{\beta}} \ar[ddd]^(.6){F_k}\\
\Min \ar[ddd]_{0} \ar@<.5ex>[rr]^(.3){\alpha\beta} && \Mout\ar@<.5ex>[ll]^(.7){\gamma} \ar[ul]_{\widetilde{\alpha}} \ar[ddd]^0
\\ \\ 
& \frac{\Ker\gamma''}{\Im\beta''}\oplus\Im\gamma''\oplus\frac{\Ker\alpha''}{\Im\gamma''}\ar[dl]_(.55){\widetilde{\beta''}} \\
\Min'' \ar@<.5ex>[rr]^{\alpha''\beta''} && \Mout''\ar@<.5ex>[ll]^{\gamma''} \ar[ul]_{\widetilde{\alpha''}}}$$
where $F_k=\left(\begin{array}{ccc} 
0 & 0 & 0\\
0 & 0 & 0 \\
\sigma\circ{\fin}\circ\varepsilon'\circ j'\circ \widetilde{{\fout'}}\circ\widetilde{\rho''} & 0 & 0 \end{array}\right).$ 
Therefore, 
the image of $\widetilde{f}\circ\widetilde{f'}-\widetilde{f\circ f'}$ lies in $\add S_k'$, 
showing that we have a functor ${\fd \Lambda} \rightarrow \frac{\fd \Lambda'}{[\add S_k']}.$

We have $F(S_k) = 0$,
hence there is a functor
$$\frac{\fd \Lambda}{[\add S_k]} \rightarrow \frac{\fd \Lambda'}{[\add S_k']}.$$

We now verify that this functor is an equivalence.  
By (the proof of) \cite[Theorem 10.13]{dwz}, this functor is a bijection at the level of representations: in particular it is dense.
To show that it is a bijection at the level of the morphisms, suppose that we have a morphism of representations 
as in diagram (\ref{commdia}).
By construction $\widetilde{\widetilde{f}}$ is of the form:
$$
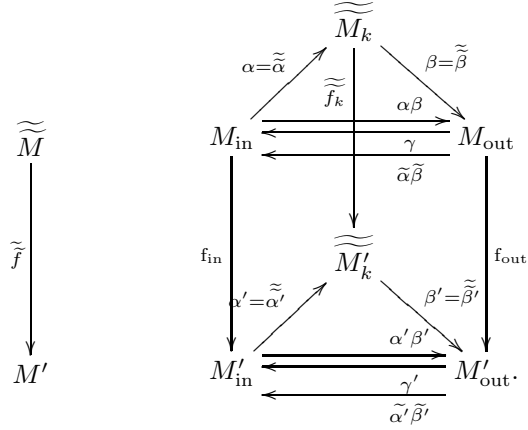
\xymatrix{&&& \widetilde{\widetilde{M_k}} \ar[dr]^{\beta=\widetilde{\widetilde{\beta}}} \ar[dd]_(.3){\widetilde{\widetilde{f_k}}}\\
\widetilde{\widetilde{M}} \ar[dd]_{\widetilde{\widetilde{f}}} && \Min \ar@<1.5ex>[rr]^(.7){\alpha\beta} 
\ar[ur]^{\alpha=\widetilde{\widetilde{\alpha}}} \ar[dd]_{\fin} 
&& \Mout \ar@<-.5ex>[ll]^(.3){\gamma} \ar@<1.5ex>[ll]^(.3){\widetilde{\alpha}\widetilde{\beta}} \ar[dd]^{\fout}
\\ &&& \widetilde{\widetilde{M'_k}} \ar[dr]^{\beta'=\widetilde{\widetilde{\beta'}}} \\
M' && \Min' \ar@<1.5ex>[rr]^(.7){\alpha'\beta'} \ar[ur]^{\alpha'=\widetilde{\widetilde{\alpha'}}} 
&& \Mout' . \ar@<-.5ex>[ll]^(.3){\gamma'} \ar@<2ex>[ll]^(.3){\widetilde{\alpha'}\widetilde{\beta'}}
}$$
By \cite{dwz}, $M$ is isomorphic to the reduced form $\mu_k(M)$ of $\widetilde{\widetilde{M}}$.  
To obtain this reduced form, one applies three steps; see \cite[Theorem 5.7, Steps 1-3]{dwz}. 
It turns out that only the first step has an impact on $\alpha, \beta$ and $\gamma$, and this modification 
consists of changing the map $\beta$ to $-\beta$. So, applying these steps does not change the 
commutativity of the above diagram, saying that after the reduction, $\widetilde{\widetilde{f}}$ is still 
a morphism for the reduced representation.

Now, since $f-\widetilde{\widetilde{f}}$ is (possibly) not zero only at position $k$, 
then $f$ and $\widetilde{\widetilde{f}}$ coincide in $\frac{\fd \Lambda}{[\add S_k]}$, showing the 
bijection at the level of morphisms.

Hence we have proved the following.

\begin{theorem}\label{nearly} 
The Jacobian algebras $\Lambda = \P(Q,W)$ and $\Lambda' = \P(\mu_k(Q',W'))$ are 
nearly Morita equivalent.
\end{theorem}

For Jacobi-finite QP's this result can be 
obtained with a very different approach by combining the result of \cite{amiot} 
saying that these algebras are 2-CY-tilted with the corresponding result for 
2-CY-tilted algebras \cite{bmr1,kr1}.

\end{document}